\documentclass[12pt]{article}%
\usepackage{amsmath}
\usepackage{amsfonts}
\usepackage{amssymb}
\usepackage{graphicx}
\usepackage{color}%
\providecommand{\U}[1]{\protect\rule{.1in}{.1in}}
\newtheorem{theorem}{Theorem}[section]

\newtheorem{corollary}[theorem]{Corollary}

\newtheorem{definition}[theorem]{Definition}

\newtheorem{lemma}[theorem]{Lemma}
\newtheorem{notation}[theorem]{Notation}

\newtheorem{remark}[theorem]{Remark}

\newenvironment{proof}[1][Proof]{\textbf{#1.} }{\ \rule{0.5em}{0.5em}}

\makeatletter\@addtoreset{equation}{section}\makeatother
\evensidemargin=\oddsidemargin
\newdimen\dummy
\dummy=\oddsidemargin
\begin{document}

\title{Critical Functions and Inf-Sup Stability of Crouzeix-Raviart Elements}
\author{C. Carstensen\thanks{(cc@math.hu-berlin.de), Department of Mathematics,
Humboldt-Universit\"{a}t zu Berlin, Berlin, D-10099, Germany.}
\and S. Sauter\thanks{(stas@math.uzh.ch), Institut f\"{u}r Mathematik,
Universit\"{a}t Z\"{u}rich, Winterthurerstr 190, CH-8057 Z\"{u}rich,
Switzerland}}
\maketitle

\begin{abstract}
In this paper, we prove that Crouzeix-Raviart finite elements of polynomial
order $p\geq5$, $p$ odd, are inf-sup stable for the Stokes problem on
triangulations. For $p\geq4$, $p$ even, the stability was proved by \'{A}.
Baran and G. Stoyan in 2007 by using the \textit{macroelement technique,} a
\textit{dimension formula}, the concept of \textit{critical points} in a
triangulation and a representation of the corresponding \textit{critical
functions}. Baran and Stoyan proved that these critical functions belong to
the range of the divergence operator applied to Crouzeix-Raviart velocity
functions and the macroelement technique implies the inf-sup stability.

The generalization of this theory to cover odd polynomial orders $p\geq5$ is
involved; one reason is that the macroelement classes, which have been used
for even $p$, are unsuitable for odd $p$. In this paper, we introduce a new
and simple representation of non-conforming Crouzeix-Raviart basis functions
of odd degree. We employ only one type of macroelement and derive
representations of all possible critical functions. Finally, we show that they
are in the range of the divergence operator applied to Crouzeix-Raviart
velocities from which the stability of the discretization follows.

\end{abstract}



\section{Introduction}

In the seminal paper \cite{CrouzeixRaviart} in 1973, Crouzeix and Raviart
developed a non-conforming family of finite elements with the goal to obtain a
stable discretization of the Stokes equation with relatively few unknowns. The
family was indexed by the (local polynomial) order, say\footnote{$\mathbb{N}%
=\left\{  1,2,\ldots\right\}  $ and $\mathbb{N}_{0}=\mathbb{N}\cup\left\{
0\right\}  $.}, $p\in\mathbb{N}$, for the approximation of the velocity field,
while the pressure is approximated by discontinuous local polynomials of
degree $p-1$. In their paper, the authors prove that for $p=1$ and spatial
dimension $d=2,3$, their non-conforming finite element leads to a stable
discretization of the Stokes equation. Since then the development of pairs of
finite elements for the velocity and for the pressure of the Stokes problem
was the topic of vivid research in numerical analysis. Surprisingly, the
problem is still not fully settled with some open issues for problems in two
spatial dimension, while much less is known for the discretization of
three-dimensional Stokes problems.

An important theoretical approach to prove the stability of Crouzeix-Raviart
elements for the Stokes problem is based on the \textit{macroelement
technique} which goes back to \cite{Stenberg_macro}, \cite{Stenberg_marco_new}%
. If the divergence operator maps a localized velocity space on the
macroelements onto the local pressure space (modulo the constant function),
the stability of the discretization follows. It has been shown in
\cite{ScottVogelius} that this surjectivity holds already for continuous
velocities for $p\geq4$ if the mesh does not contain \textit{critical points}
(cf. Def. \ref{DefCritpoint}). If the mesh contains critical points an
established idea to prove inf-sup stability is as follows: first, one
identifies the \textit{critical functions} whose span has zero intersection
with the range of the local divergence operator applied to continuous
velocities. Then, one proves that the critical functions lie in the span of
the divergence operator applied to Crouzeix-Raviart velocities and stability
of the discretization follows.

The case for even $p\geq4$ was proved along these lines in \cite{Baran_Stoyan}%
. In that paper, seven types of non-overlapping macroelements are considered
consisting of two and three triangles, the critical functions are identified,
and it is shown that these belong to the range of the local divergence
operator applied to the localized Crouzeix-Raviart velocities so that inf-sup
stability follows for this case. In our paper, we focus on the case of odd
$p\geq5$ and proceed conceptually in the same way. However, it turns out that
the macroelements, which are considered in \cite{Baran_Stoyan}, are not suited
for odd $p$. One reason is that the non-conforming basis functions for even
$p$ are more local (support on one triangle) compared to odd $p$ (support on
two adjacent triangles). Instead, we consider here nodal patches (element
stars) for the interior vertices of the triangulation as the only type of
overlapping macroelements. We identify the \textit{critical functions} in
\S \ref{PTpge5}, which are related to \textit{critical points} (cf. Def.
\ref{DefCritpoint}) in the nodal patches. We show that these critical
functions form a basis of the complement in the pressure space of the range of
the divergence operator applied to the localized continuous velocities; the
proof employs the \textquotedblleft dimension formula\textquotedblright\ in
\cite{ScottVogelius} and hence is restricted to $p\geq4$. Then we show that
the critical functions belong to the range of the divergence operator applied
to Crouzeix-Raviart velocities and this implies the stability of the discretization.

The main achievements in this paper are as follows. A new and simple
representation of Crouzeix-Raviart basis functions for odd $p$ is introduced
in Section \ref{VelBasis}. We identify the critical pressure functions for odd
$p\geq5$ and continuous localized velocities in Definition \ref{Defcritfunc}
and show that they form a basis for the complement of the range of the
divergence operator applied to localized continuous velocities. Finally, we
prove that the critical functions belong to the range of the divergence
operator applied to localized Crouzeix-Raviart velocities.

This leads to the main conclusion: for $p\in\mathbb{N}$, let $CR_{p,0}\left(
\mathcal{T}\right)  $ denote the (scalar) Crouzeix-Raviart finite element
space with local polynomial order $p$ on regular triangulations $\mathcal{T}$,
i.e., without hanging nodes, of a two-dimensional bounded polygonal Lipschitz
domain $\Omega\subset\mathbb{R}^{2}$ obeying zero-boundary conditions in a
\textquotedblleft Crouzeix-Raviart\textquotedblright\ sense. Let
$\mathbb{P}_{p-1}\left(  \mathcal{T}\right)  $ denote the discontinuous finite
element space of local polynomial order $p-1$ on this triangulation. Then%
\begin{equation}
\left(  \left(  CR_{p,0}\left(  \mathcal{T}\right)  \right)  ^{2}%
,\mathbb{P}_{p-1}\left(  \mathcal{T}\right)  /\mathbb{R}\right)  \text{ is a
stable finite element for the Stokes equation} \label{mainresult}%
\end{equation}
for odd $p\geq5$.

In the following we comment on the remaining cases: $p=1,2,3$ and even
$p\geq4$. In \cite{CrouzeixRaviart}, statement (\ref{mainresult}) is proved
for $p=1$ (and for spatial dimensions $d=2,3$). From \cite{Baran_Stoyan} we
know that assertion (\ref{mainresult}) is true for even $p\geq4$. For $p=2$,
the result we proved already in \cite{Fortin_Soulie}. For the case $p\geq5$
and $p$ odd, statement (\ref{mainresult}) follows from \cite{ScottVogelius}
(see also \cite{GuzmanScott2019}) provided the triangulation does not contain
critical points. In \cite{Crouzeix_Falk}, the case $p=3$ is considered and the
claim (\ref{mainresult}) is proved if the triangulation does not contain
critical points\ and an additional technical condition for the nodal points is
satisfied. The case $p=3$ of statement (\ref{mainresult}) for any regular
triangulation is proved in \cite{CCSS_CR_2} which is the companion to this
paper. The proof in \cite{CCSS_CR_2} circumvents the concept of critical
points and functions and reduces the problem to a purely algebraic problem in
that to determe the nullspace for a coefficient matrix of a linear system on
the local nodal patches. The proof in \cite{CCSS_CR_2} applies also to odd
$p\geq3$ but does not give insight on the mechanism how the critical functions
are eliminated by the non-conforming Crouzeix-Raviart functions.

The paper is structured as follows. In Section \ref{SecFEM} we formulate the
Stokes equation in variational form and introduce the functional analytic
setting in \S \ref{SecContStokes}. The discretization is based on
non-conforming Galerkin finite element methods on regular triangulations. The
finite element spaces, in particular the Crouzeix-Raviart finite elements of
polynomial order $p$ are introduced in Section \ref{NuDiscrete}. The inf-sup
condition and the main theorem (Thm. \ref{Theomain}) are stated at the end of
this section.

The final section \S \ref{SecProof} of the paper is devoted to the proof of
the main theorem. In the appendices \S \ref{DerBary}, \S \ref{SecMatDet},
\S \ref{SecIntJP} we provide some technical properties on derivatives of
barycentric coordinates, determinants of tridiagonal matrices, and closed form
integrals of some products of Jacobi polynomials.

\section{Setting\label{SecFEM}}

\subsection{The Continuous Stokes Problem\label{SecContStokes}}

Let $\Omega\subset\mathbb{R}^{2}$ denote a bounded polygonal domain with
boundary $\partial\Omega$. Our goal is to find a family of pairs of finite
element spaces for the stable numerical solution of the Stokes equation. On
the continuous level, the strong form of the Stokes equation is given by%
\[%
\begin{array}
[c]{llll}%
-\Delta\mathbf{u} & +\nabla p & =f & \text{in }\Omega,\\
\operatorname*{div}\mathbf{u} &  & =0 & \text{in }\Omega
\end{array}
\]
with boundary conditions for the velocity and a normalization condition for
the pressure%
\[
\mathbf{u}=\mathbf{0}\quad\text{on }\partial\Omega\quad\text{and\quad}%
\int_{\Omega}p=0.
\]
To state the classical existence and uniqueness result we formulate this
equation in a variational form and first introduce the relevant function
spaces. Throughout the paper we restrict to vector spaces over the field of
real numbers.

For $s\geq0$, $1\leq p\leq\infty$, $W^{s,p}\left(  \Omega\right)  $ denote the
classical Sobolev spaces of functions with norm $\left\Vert \cdot\right\Vert
_{W^{s,p}\left(  \Omega\right)  }$. As usual we write $L^{p}\left(
\Omega\right)  $ instead of $W^{0,p}\left(  \Omega\right)  $ and $H^{s}\left(
\Omega\right)  $ for $W^{s,2}\left(  \Omega\right)  $. For $s\geq0$, we denote
by $H_{0}^{s}\left(  \Omega\right)  $ the closure with respect to the
$H^{s}\left(  \Omega\right)  $ norm of the space of infinitely smooth
functions with compact support in $\Omega$. Its dual space is denoted by
$H^{-s}\left(  \Omega\right)  $.

The scalar product and norm in $L^{2}\left(  \Omega\right)  $ are denoted
respectively by
%\ednote{Ivan check if we need the inner products}%
\[%
\begin{array}
[c]{llll}%
\left(  u,v\right)  _{L^{2}\left(  \Omega\right)  }:=\int_{\Omega}uv &
\text{and} & \left\Vert u\right\Vert _{L^{2}\left(  \Omega\right)  }:=\left(
u,u\right)  _{L^{2}\left(  \Omega\right)  }^{1/2} & \text{in }L^{2}\left(
\Omega\right)  .
\end{array}
\]
Vector-valued and $2\times2$ tensor-valued analogues of the function spaces
are denoted by bold and blackboard bold letters, e.g., $\mathbf{H}^{s}\left(
\Omega\right)  =\left(  H^{s}\left(  \Omega\right)  \right)  ^{2}$ and
$\mathbb{H}^{s}=\left(  H^{s}\left(  \Omega\right)  \right)  ^{2\times2}$ and
analogously for other quantities.

The $\mathbf{L}^{2}\left(  \Omega\right)  $ scalar product and norm for vector
valued functions are given by%
\[
\left(  \mathbf{u},\mathbf{v}\right)  _{\mathbf{L}^{2}\left(  \Omega\right)
}:=\int_{\Omega}\left\langle \mathbf{u},\mathbf{v}\right\rangle \quad
\text{and\quad}\left\Vert \mathbf{u}\right\Vert _{\mathbf{L}^{2}\left(
\Omega\right)  }:=\left(  \mathbf{u},\mathbf{u}\right)  _{\mathbf{L}%
^{2}\left(  \Omega\right)  }^{1/2},
\]
where $\left\langle \mathbf{u},\mathbf{v}\right\rangle $ denotes the Euclidean
scalar product in $\mathbb{R}^{2}$. In a similar fashion, we define for
$\mathbf{G},\mathbf{H}\in\mathbb{L}^{2\times2}\left(  \Omega\right)  $ the
scalar product and norm by%
\[
\left(  \mathbf{G},\mathbf{H}\right)  _{\mathbb{L}^{2\times2}\left(
\Omega\right)  }:=\int_{\Omega}\left\langle \mathbf{G},\mathbf{H}\right\rangle
\quad\text{and\quad}\left\Vert \mathbf{G}\right\Vert _{\mathbb{L}^{2\times
2}\left(  \Omega\right)  }:=\left(  \mathbf{G},\mathbf{G}\right)
_{\mathbb{L}^{2\times2}\left(  \Omega\right)  }^{1/2},
\]
where $\left\langle \mathbf{G},\mathbf{H}\right\rangle =\sum_{i,j=1}%
^{2}G_{i,j}H_{i,j}$. Finally, let $L_{0}^{2}\left(  \Omega\right)  :=\left\{
u\in L^{2}\left(  \Omega\right)  :\int_{\Omega}u=0\right\}  $.

We introduce the bilinear form $a:\mathbf{H}^{1}\left(  \Omega\right)
\times\mathbf{H}^{1}\left(  \Omega\right)  \rightarrow\mathbb{R}$ by%
\[
a\left(  \mathbf{u},\mathbf{v}\right)  :=\left(  \nabla\mathbf{u}%
,\nabla\mathbf{v}\right)  _{\mathbb{L}^{2\times2}\left(  \Omega\right)  },
\]
where $\nabla\mathbf{u}$ and $\nabla\mathbf{v}$ denote the derivatives (Jacobi
matrices) of $\mathbf{u}$ and $\mathbf{v}$. The variational form of the Stokes
problem is given by: For given $\mathbf{f}\in\mathbf{H}^{-1}\left(
\Omega\right)  ,$%
\begin{equation}
\text{find }\left(  \mathbf{u},p\right)  \in\mathbf{H}_{0}^{1}\left(
\Omega\right)  \times L_{0}^{2}\left(  \Omega\right)  \;\text{s.t.\ }\left\{
\begin{array}
[c]{lll}%
a\left(  \mathbf{u},\mathbf{v}\right)  -\left(  p,\operatorname*{div}%
\mathbf{v}\right)  _{L^{2}\left(  \Omega\right)  } & =\left(  \mathbf{f}%
,\mathbf{v}\right)  _{\mathbf{L}^{2}\left(  \Omega\right)  } & \forall
\mathbf{v}\in\mathbf{H}_{0}^{1}\left(  \Omega\right)  ,\\
\left(  \operatorname*{div}\mathbf{u},q\right)  _{L^{2}\left(  \Omega\right)
} & =0 & \forall q\in L_{0}^{2}\left(  \Omega\right)  .
\end{array}
\right.  \label{varproblemstokes}%
\end{equation}

It is well-known (see, e.g., \cite{Girault86}) that (\ref{varproblemstokes})
is well posed.

\subsection{Numerical Discretization of the Stokes Problem\label{NuDiscrete}}

In the following we introduce a discretization for problem
(\ref{varproblemstokes}). Let $\mathcal{T}=\left\{  K_{i}:1\leq i\leq
n\right\}  $ denote a \textit{regular }triangulation of $\Omega$ consisting of
closed triangles $K_{i}$ which have the property that the intersection of two
different triangles $K_{i}$, $K_{j}$ is either empty, a common edge, or a
common point. We also assume $\Omega=\operatorname*{dom}\mathcal{T}$, where
\begin{equation}
\operatorname*{dom}\mathcal{T}:=\operatorname*{int}\left(
%TCIMACRO{\dbigcup \limits_{K\in\mathcal{T}}}%
%BeginExpansion
{\displaystyle\bigcup\limits_{K\in\mathcal{T}}}
%EndExpansion
K\right)  \label{defdomT}%
\end{equation}
and $\operatorname*{int}\left(  M\right)  $ denotes the interior of a set
$M\subset\mathbb{R}^{2}$. An important measure for the quality of a finite
element triangulation is the shape-regularity constant, which we define by%
\begin{equation}
\gamma_{\mathcal{T}}:=\max_{K\in\mathcal{T}}\frac{h_{K}}{\rho_{K}}
\label{defgammat}%
\end{equation}
with the local mesh width $h_{K}:=\operatorname*{diam}K$ and $\rho_{K}$
denoting the diameter of the largest inscribed ball in $K$.

The set of edges in $\mathcal{T}$ are denoted by $\mathcal{E}$, while the
subset of boundary edges is $\mathcal{E}_{\partial\Omega}:=\left\{
E\in\mathcal{E}:E\subset\partial\Omega\right\}  $; the subset of inner edges
is given by $\mathcal{E}_{\Omega}:=\mathcal{E}\backslash\mathcal{E}%
_{\partial\Omega}$. The set of triangle vertices in $\mathcal{T}$ is denoted
by $\mathcal{V}$, while the subset of inner vertices is $\mathcal{V}_{\Omega
}:=\left\{  \mathbf{V}\in\mathcal{V}:\mathbf{V}\notin\partial\Omega\right\}  $
and $\mathcal{V}_{\partial\Omega}:=\mathcal{V}\backslash\mathcal{V}_{\Omega}$.
For $E\in\mathcal{E}$, we define the edge patch by%
\[
\mathcal{T}_{E}:=\left\{  K\in\mathcal{T}:E\subset K\right\}  \quad
\text{and\quad}\omega_{E}:=%
%TCIMACRO{\dbigcup \limits_{K\in\mathcal{T}_{E}}}%
%BeginExpansion
{\displaystyle\bigcup\limits_{K\in\mathcal{T}_{E}}}
%EndExpansion
K.
\]
For $\mathbf{z}\in\mathcal{V}$, the nodal patch is defined by%
\begin{equation}
\mathcal{T}_{\mathbf{z}}:=\left\{  K\in\mathcal{T}:\mathbf{z}\in K\right\}
\quad\text{and\quad}\omega_{\mathbf{z}}:=%
%TCIMACRO{\dbigcup \limits_{K\in\mathcal{T}_{\mathbf{z}}}}%
%BeginExpansion
{\displaystyle\bigcup\limits_{K\in\mathcal{T}_{\mathbf{z}}}}
%EndExpansion
K. \label{nodalpatch}%
\end{equation}

We will need an additional mesh parameter. For a regular triangulation
$\mathcal{T}$ of $\Omega$, let%
\[
\mathcal{T}^{\prime}:=%
%TCIMACRO{\dbigcup \limits_{\mathbf{z}\in\mathcal{V}_{\Omega}}}%
%BeginExpansion
{\displaystyle\bigcup\limits_{\mathbf{z}\in\mathcal{V}_{\Omega}}}
%EndExpansion
\mathcal{T}_{\mathbf{z}}.
\]
Then, the constant\footnote{By $\left\vert \mathcal{J}\right\vert $ we denote
the cardinality of a discrete set $\mathcal{J}$ (cf. Notation \ref{Notation}%
).}%
\begin{equation}
d_{\mathcal{T}}:=\left\vert \mathcal{T}\backslash\mathcal{T}^{\prime
}\right\vert \label{dtau}%
\end{equation}
denotes the number of triangles in the triangulation which are not connected
to an inner point.

For $m\in\mathbb{N}$, we employ the usual multiindex notation for $%
%TCIMACRO{\TeXButton{boldmu}{\mbox{\boldmath$ \mu$}}}%
%BeginExpansion
\mbox{\boldmath$ \mu$}%
%EndExpansion
=\left(  \mu_{i}\right)  _{i=1}^{m}\in\mathbb{N}_{0}^{m}$ and points
$\mathbf{x}=\left(  x_{i}\right)  _{i=1}^{m}\in\mathbb{R}^{m}$%
\[
\left\vert
%TCIMACRO{\TeXButton{boldmu}{\mbox{\boldmath$ \mu$}}}%
%BeginExpansion
\mbox{\boldmath$ \mu$}%
%EndExpansion
\right\vert :=\mu_{1}+\ldots+\mu_{m},\quad\mathbf{x}^{%
%TCIMACRO{\TeXButton{boldmu}{\mbox{\boldmath$ \mu$}}}%
%BeginExpansion
\mbox{\boldmath$ \mu$}%
%EndExpansion
}:=%
%TCIMACRO{\dprod \limits_{j=1}^{m}}%
%BeginExpansion
{\displaystyle\prod\limits_{j=1}^{m}}
%EndExpansion
x_{j}^{\mu_{j}}.
\]

For $k\in\mathbb{N}_{0}$ and $m\in\mathbb{N}$, we define the index sets%
\[
\mathbb{I}_{\leq k}^{m}:=\left\{
%TCIMACRO{\TeXButton{boldmu}{\mbox{\boldmath$ \mu$}}}%
%BeginExpansion
\mbox{\boldmath$ \mu$}%
%EndExpansion
\in\mathbb{N}_{0}^{m}\mid\left\vert
%TCIMACRO{\TeXButton{boldmu}{\mbox{\boldmath$ \mu$}}}%
%BeginExpansion
\mbox{\boldmath$ \mu$}%
%EndExpansion
\right\vert \leq k\right\}  \quad\text{and\quad}\mathbb{I}_{=k}^{m}:=\left\{
%TCIMACRO{\TeXButton{boldmu}{\mbox{\boldmath$ \mu$}}}%
%BeginExpansion
\mbox{\boldmath$ \mu$}%
%EndExpansion
\in\mathbb{N}_{0}^{m}\mid\left\vert
%TCIMACRO{\TeXButton{boldmu}{\mbox{\boldmath$ \mu$}}}%
%BeginExpansion
\mbox{\boldmath$ \mu$}%
%EndExpansion
\right\vert =k\right\}
\]
Let $\mathbb{P}_{m,k}$ denote the space of $m$-variate polynomials of maximal
degree $k$, consisting of functions of the form%
\[
\sum_{%
%TCIMACRO{\TeXButton{boldmu}{\mbox{\boldmath$ \mu$}}}%
%BeginExpansion
\mbox{\boldmath$ \mu$}%
%EndExpansion
\in\mathbb{I}_{\leq k}^{m}}a_{%
%TCIMACRO{\TeXButton{boldmu}{\mbox{\boldmath$ \mu$}}}%
%BeginExpansion
\mbox{\boldmath$ \mu$}%
%EndExpansion
}\mathbf{x}^{%
%TCIMACRO{\TeXButton{boldmu}{\mbox{\boldmath$ \mu$}}}%
%BeginExpansion
\mbox{\boldmath$ \mu$}%
%EndExpansion
}%
\]
for real coefficients $a_{%
%TCIMACRO{\TeXButton{boldmu}{\mbox{\boldmath$ \mu$}}}%
%BeginExpansion
\mbox{\boldmath$ \mu$}%
%EndExpansion
}$. Formally, we set $\mathbb{P}_{m,-1}:=\left\{  0\right\}  $. To indicate
the domain of definition we write sometimes $\mathbb{P}_{k}\left(  D\right)  $
for $D\subset\mathbb{R}^{m}$ and skip the index $m$ since it is then clear
from the argument $D$.

For $s\geq0$ and a regular triangulation $\mathcal{T}$ for the domain $\Omega
$, let%
\[
H^{s}\left(  \mathcal{T}\right)  :=\left\{  u\in L^{2}\left(  \Omega\right)
\mid\forall K\in\mathcal{T}:\left.  u\right\vert _{K}\in H^{1}\left(
K\right)  \right\}  .
\]
We introduce the following finite element spaces%
\begin{equation}%
\begin{array}
[c]{ll}
& \mathbb{P}_{k}\left(  \mathcal{T}\right)  :=\left\{  q\in L^{2}\left(
\Omega\right)  \mid\forall K\in\mathcal{T}:\left.  q\right\vert _{K}%
\in\mathbb{P}_{k}\left(  K\right)  \right\}  ,\\
\text{and (cf. (\ref{defdomT}))} & \mathbb{P}_{k}\left(  \mathcal{T}\right)
/\mathbb{R}:=\left\{  q\in\mathbb{P}_{k}\left(  \mathcal{T}\right)
:\int_{\operatorname*{dom}\mathcal{T}}q=0\right\}  .
\end{array}
\label{Pkdefs}%
\end{equation}
Furthermore, let%
\[%
\begin{array}
[c]{ll}
& S_{k}\left(  \mathcal{T}\right)  :=\left\{  v\in C^{0}\left(  \Omega\right)
\mid\forall K\in\mathcal{T}:\left.  v\right\vert _{K}\in\mathbb{P}_{k}\left(
K\right)  \right\}  ,\\
\text{and} & S_{k,0}\left(  \mathcal{T}\right)  :=S_{k}\left(  \mathcal{T}%
\right)  \cap H_{0}^{1}\left(  \operatorname*{dom}\mathcal{T}\right)  .
\end{array}
\]
The vector-valued versions are denoted by $\mathbf{S}_{k}\left(
\mathcal{T}\right)  :=S_{k}\left(  \mathcal{T}\right)  ^{2}$ and
$\mathbf{S}_{k,0}\left(  \mathcal{T}\right)  :=S_{k,0}\left(  \mathcal{T}%
\right)  ^{2}$. Finally, we define the Crouzeix-Raviart space by%
%TCIMACRO{\TeXButton{CRdeffull}{\begin{subequations}
%\label{CRdeffull}
%\end{subequations}}}%
%BeginExpansion
\begin{subequations}
\label{CRdeffull}
\end{subequations}%
%EndExpansion%
\begin{align}
CR_{k}\left(  \mathcal{T}\right)   &  :=\left\{  v\in\mathbb{P}_{k}\left(
\mathcal{T}\right)  \mid\forall q\in\mathbb{P}_{k-1}\left(  E\right)
\quad\forall E\in\mathcal{E}_{\Omega}\quad\int_{E}\left[  v\right]
_{E}q=0\right\}  ,\tag{%
%TCIMACRO{\TeXButton{CRdeffull}{\ref{CRdeffull}}}%
%BeginExpansion
\ref{CRdeffull}%
%EndExpansion
a}\label{CRdeffulla}\\
CR_{k,0}\left(  \mathcal{T}\right)   &  :=\left\{  v\in CR_{k}\left(
\mathcal{T}\right)  \mid\forall q\in\mathbb{P}_{k-1}\left(  E\right)
\quad\forall E\in\mathcal{E}_{\partial\Omega}\quad\int_{E}vq=0\right\}  .
\tag{%
%TCIMACRO{\TeXButton{CRdeffull}{\ref{CRdeffull}}}%
%BeginExpansion
\ref{CRdeffull}%
%EndExpansion
b}\label{PCR0RB}%
\end{align}
Here, $\left[  v\right]  _{E}$ denotes the jump of $v\in\mathbb{P}_{k}\left(
\mathcal{T}\right)  $ across an edge $E\in\mathcal{E}_{\Omega}$ and
$\mathbb{P}_{k-1}\left(  E\right)  $ is the space of polynomials of maximal
degree $k-1$ with respect to the local variable in $E$.

We have collected all ingredients for defining the Crouzeix-Raviart
discretization for the Stokes equation. For $p\in\mathbb{N}$, let the discrete
velocity space and pressure space be defined by
\[
\mathbf{CR}_{p,0}\left(  \mathcal{T}\right)  :=\left(  CR_{p,0}\left(
\mathcal{T}\right)  \right)  ^{2}\quad\text{and\quad}M_{p-1}\left(
\mathcal{T}\right)  :=\mathbb{P}_{p-1}\left(  \mathcal{T}\right)
/\mathbb{R}.
\]
Then, the discretization is given by: find $\left(  \mathbf{u}%
_{\operatorname*{CR}},p_{\operatorname*{disc}}\right)  \in\mathbf{CR}%
_{p,0}\left(  \mathcal{T}\right)  \times M_{p-1}\left(  \mathcal{T}\right)
\;$such that%
\begin{equation}
\left\{
\begin{array}
[c]{lll}%
a\left(  \mathbf{u}_{\operatorname*{CR}},\mathbf{v}\right)  -\left(
p_{\operatorname*{disc}},\operatorname*{div}\mathbf{v}\right)  _{L^{2}\left(
\Omega\right)  } & =\left(  \mathbf{f},\mathbf{v}\right)  _{\mathbf{L}%
^{2}\left(  \Omega\right)  } & \forall\mathbf{v}\in\mathbf{CR}_{p,0}\left(
\mathcal{T}\right)  ,\\
\left(  \operatorname*{div}\mathbf{u}_{\operatorname*{CR}},q\right)
_{L^{2}\left(  \Omega\right)  } & =0 & \forall q\in M_{p-1}\left(
\mathcal{T}\right)  .
\end{array}
\right.  \label{discrStokes}%
\end{equation}

\begin{definition}
Let $\mathcal{T}$ denote a regular triangulation for $\Omega$. A pair
$\mathbf{CR}_{p,0}\left(  \mathcal{T}\right)  \times M_{p-1}\left(
\mathcal{T}\right)  $ is \emph{inf-sup stable} if there exists a constant
$c_{\mathcal{T},p}$ such that%
\begin{equation}
\inf_{p\in M_{p-1}\left(  \mathcal{T}\right)  \backslash\left\{  0\right\}
}\sup_{\mathbf{v}\in\mathbf{S}_{p,0}^{\operatorname*{CR}}\left(
\mathcal{T}\right)  \backslash\left\{  0\right\}  }\frac{\left(
p,\operatorname*{div}\mathbf{v}\right)  _{L^{2}\left(  \Omega\right)  }%
}{\left\Vert \mathbf{v}\right\Vert _{\mathbf{H}^{1}\left(  \Omega\right)
}\left\Vert p\right\Vert _{L^{2}\left(  \Omega\right)  }}\geq c_{\mathcal{T}%
,p}>0. \label{infsupcond}%
\end{equation}

\end{definition}

We are now in the position to formulate our main theorem.

\begin{theorem}
\label{Theomain}Let $\Omega\subset\mathbb{R}^{2}$ be a bounded polygonal
Lipschitz domain and let $\mathcal{T}$ denote a regular triangulation of
$\Omega$, which contains at least one inner node. Then, the inf-sup condition
(\ref{infsupcond}) holds for a constant $c_{\mathcal{T},p}$, which depends on
the shape regularity of the mesh, the constant $d_{\mathcal{T}}$ in
(\ref{dtau}), and the polynomial degree $p$.
\end{theorem}

We emphasize that the original definition in \cite{CrouzeixRaviart} allows for
slightly more general finite element spaces, more precisely, the spaces
$CR_{p}\left(  \mathcal{T}\right)  $ can be enriched by locally supported
functions. From this point of view, the definition (\ref{CRdeffull}) describes
the \textit{minimal} Crouzeix-Raviart space. The possibility for enrichment
has been used frequently in the literature to prove inf-sup stability for the
arising finite element spaces (see, e.g., \cite{CrouzeixRaviart},
\cite{Guzman_divfree}, \cite{Matthies_nonconf_2005}). In contrast, we will
prove the stability for the minimal Crouzeix-Raviart family.

\section{Proof of Theorem \ref{Theomain}\label{SecProof}}

In \cite{CrouzeixRaviart}, Theorem \ref{Theomain} is proved for $p=1$ (and for
spatial dimensions $d=2,3$). From \cite{Baran_Stoyan} we know that the theorem
is true for even $p\geq4$. In \cite{Fortin_Soulie}, the result is proved for
$p=2$. In this section, we will prove the result for odd $p\geq5$ and refer
for the proof of the case $p=3$ to \cite{CCSS_CR_2}.

\subsection{Barycentric Coordinates and Basis Functions for the
Velocity\label{VelBasis}}

In this section, we introduce basis functions for the finite element spaces in
Section \ref{NuDiscrete}. We begin with introducing some general notation.

\begin{notation}
\label{Notation}For vectors $\mathbf{a}_{i}\in\mathbb{R}^{n}$, $1\leq i\leq
m$, we write $\left[  \mathbf{a}_{1}\mid\mathbf{a}_{2}\mid\ldots\mid
\mathbf{a}_{m}\right]  $ for the $n\times m$ matrix with column vectors
$\mathbf{a}_{i}$. For $\mathbf{v}=\binom{v_{1}}{v_{2}}\in\mathbb{R}^{2}$ we
set $\mathbf{v}^{\perp}=\left(  v_{2},-v_{1}\right)  ^{T}$. Let\ $\mathbf{e}%
_{k,i}\in\mathbb{R}^{k}$\ be the $i$-th canonical unit vector in
$\mathbb{R}^{k}.$

Vertices in a triangle are always numbered counterclockwise. In a triangle $K$
with vertices $\mathbf{A}_{1}$, $\mathbf{A}_{2}$, $\mathbf{A}_{3}$ the angle
at $\mathbf{A}_{i}$ is called $\alpha_{i}$ or alternatively $\alpha_{j,k}$
where $i,j,k\in\left\{  1,2,3\right\}  $ are pairwise different. If a triangle
is numbered by an index (e.g., $K_{\ell}$), the angle at $A_{\ell,i}$ is
called $\alpha_{\ell,i}$ or alternatively $\alpha_{\ell,j,k}$. For quantities
in a triangle $K$ as, e.g., angles $\alpha_{j}$, $1\leq j\leq3$, we use the
cyclic numbering convention $\alpha_{3+1}:=\alpha_{1}$ and $\alpha
_{1-1}:=\alpha_{3}$.

For a $d$-dimensional measurable set $D,$ we write $\left\vert D\right\vert $
for its measure; for a discrete set, say $\mathcal{J}$, we denote by
$\left\vert \mathcal{J}\right\vert $ its cardinality.

In the proofs, we consider frequently nodal patches $\mathcal{T}_{\mathbf{z}}$
for inner vertices $\mathbf{z}\in\mathcal{V}_{\Omega}$. The number $m$ denotes
the number of triangles in $\mathcal{T}_{\mathbf{z}}$. Various quantities in
this patch such as, e.g., the triangles in $\mathcal{T}_{\mathbf{z}}$, have an
index which runs from $1$ to $m$. Here, we use the cyclic numbering convention
$K_{m+1}:=K_{1}$ and $K_{1-1}:=K_{m}$ and apply this analogously for other
quantities in the nodal patch.
\end{notation}

Let the closed reference triangle $\widehat{K}$ be the triangle with vertices
$\mathbf{\hat{A}}_{1}:=\left(  0,0\right)  $, $\mathbf{\hat{A}}_{2}:=\left(
1,0\right)  $, $\mathbf{\hat{A}}_{3}:=\left(  0,1\right)  $. The nodal points
on the reference element of order $k\in\mathbb{N}_{0}$ are given by%
\[
\widehat{\mathcal{N}}_{k}:=\left\{
\begin{array}
[c]{ll}%
\left\{  \dfrac{1}{k}%
%TCIMACRO{\TeXButton{boldmu}{\mbox{\boldmath$ \mu$}}}%
%BeginExpansion
\mbox{\boldmath$ \mu$}%
%EndExpansion
\mid%
%TCIMACRO{\TeXButton{boldmu}{\mbox{\boldmath$ \mu$}}}%
%BeginExpansion
\mbox{\boldmath$ \mu$}%
%EndExpansion
\in\mathbb{I}_{\leq k}^{2}\right\}  & k\geq1,\\
& \\
\left\{  \left(  \dfrac{1}{3},\dfrac{1}{3}\right)  \right\}  & k=0.
\end{array}
\right.
\]
For a triangle $K\subset\mathbb{R}^{2}$, we denote by $\chi_{K}:\widehat{K}%
\rightarrow K$ an affine bijection. The mapped nodal points of order
$k\in\mathbb{N}_{0}$ on $K$ are given by%
\[
\mathcal{N}_{k}\left(  K\right)  :=\left\{  \chi_{K}\left(  \mathbf{z}\right)
:\mathbf{z}\in\widehat{\mathcal{N}}_{k}\right\}  .
\]
The nodal points of order $k$ on $\mathcal{T}$ are defined by%
\[
\mathcal{N}_{k}\left(  \mathcal{T}\right)  :=%
%TCIMACRO{\dbigcup \limits_{K\in\mathcal{T}}}%
%BeginExpansion
{\displaystyle\bigcup\limits_{K\in\mathcal{T}}}
%EndExpansion
\mathcal{N}_{k}\left(  K\right)  \quad\text{and\quad}\mathcal{N}%
_{\partial\Omega}^{k}\left(  \mathcal{T}\right)  :=\mathcal{N}_{k}\left(
\mathcal{T}\right)  \cap\partial\Omega,\quad\mathcal{N}_{k,\Omega}\left(
\mathcal{T}\right)  :=\mathcal{N}_{k}\left(  \mathcal{T}\right)  \cap\Omega.
\]
We introduce the well-known Lagrange basis for the space $S_{k}\left(
\mathcal{T}\right)  $, which is indexed by the nodal points $\mathbf{z}%
\in\mathcal{N}_{k}\left(  \mathcal{T}\right)  $ and characterized by
\begin{equation}
B_{k,\mathbf{z}}\in S_{k}\left(  \mathcal{T}\right)  \quad\text{and\quad
}\forall\mathbf{z}^{\prime}\in\mathcal{N}_{k}\left(  \mathcal{T}\right)
\qquad B_{k,\mathbf{z}}\left(  \mathbf{z}^{\prime}\right)  =\delta
_{\mathbf{z},\mathbf{z}^{\prime}}, \label{basisfunctions}%
\end{equation}
where $\delta_{\mathbf{z},\mathbf{z}^{\prime}}$ is the Kronecker delta. A
basis for the space $S_{k,0}\left(  \mathcal{T}\right)  $ is given by
$B_{k,\mathbf{z}}$, $\mathbf{z}\in\mathcal{N}_{k,\Omega}\left(  \mathcal{T}%
\right)  $.\bigskip

Next, we define a basis for the Crouzeix-Raviart space. Let $\alpha,\beta>-1$
and $n\in\mathbb{N}_{0}$. The \emph{Jacobi polynomial} $P_{n}^{\left(
\alpha,\beta\right)  }$ is a polynomial of degree $n$ such that
\[
\int_{-1}^{1}P_{n}^{\left(  \alpha,\beta\right)  }\left(  x\right)  \,q\left(
x\right)  \left(  1-x\right)  ^{\alpha}\left(  1+x\right)  ^{\beta}\,dx=0
\]
for all polynomials $q$ of degree less than $n$, and (cf. \cite[Table
18.6.1]{NIST:DLMF})%
\begin{equation}
P_{n}^{\left(  \alpha,\beta\right)  }\left(  1\right)  =\frac{\left(
\alpha+1\right)  _{n}}{n!},\qquad P_{n}^{\left(  \alpha,\beta\right)  }\left(
-1\right)  =\left(  -1\right)  ^{n}\frac{\left(  \beta+1\right)  _{n}}{n!}.
\label{Pnormalization}%
\end{equation}
Here the \emph{shifted factorial} is defined by $\left(  a\right)
_{n}:=a\left(  a+1\right)  \ldots\left(  a+n-1\right)  $ for $n>0$ and
$\left(  a\right)  _{0}:=1$. Note that $P_{k}^{\left(  0,0\right)  }$ are the
Legendre polynomials (see \cite[18.7.9]{NIST:DLMF}).\medskip

Let $K$ denote a triangle with vertices $\mathbf{A}_{i}$, $1\leq i\leq3$, and
let $\lambda_{K,\mathbf{A}_{i}}\in\mathbb{P}_{1}\left(  K\right)  $ be the
\textit{barycentric coordinate} for the node $\mathbf{A}_{i}$ defined by%
\begin{equation}
\lambda_{K,\mathbf{A}_{i}}\left(  \mathbf{A}_{j}\right)  =\delta_{i,j}%
\quad1\leq i,j\leq3. \label{lambdaintro1}%
\end{equation}
If the numbering of the vertices in $K$ is fixed, we write $\lambda_{K,i}$
short for $\lambda_{K,\mathbf{A}_{i}}$ and for $%
%TCIMACRO{\TeXButton{boldmu}{\mbox{\boldmath$ \mu$}}}%
%BeginExpansion
\mbox{\boldmath$ \mu$}%
%EndExpansion
\in\mathbb{N}_{0}^{3}$:%
\begin{equation}%
%TCIMACRO{\TeXButton{boldlambda}{\mbox{\boldmath$ \lambda$}}}%
%BeginExpansion
\mbox{\boldmath$ \lambda$}%
%EndExpansion
_{K}^{%
%TCIMACRO{\TeXButton{boldmu}{\mbox{\boldmath$ \mu$}}}%
%BeginExpansion
\mbox{\boldmath$ \mu$}%
%EndExpansion
}=\lambda_{K,1}^{\mu_{1}}\lambda_{K,2}^{\mu_{2}}\lambda_{K,3}^{\mu_{3}}.
\label{lambdaintro2}%
\end{equation}
For the barycentric coordinate on the reference element $\widehat{K}$ for the
vertex $\mathbf{\hat{A}}_{j}$ we write $\widehat{\lambda}_{j}$, $j=1,2,3$.

\begin{definition}
Let $p\in\mathbb{N}$ be even and $K\in\mathcal{T}$. Then, the
\emph{non-conforming triangle bubble }is given by%
\[
B_{p,K}^{\operatorname*{CR}}:=\left\{
\begin{array}
[c]{ll}%
\dfrac{1}{2}\left(  -1+%
%TCIMACRO{\dsum \limits_{i=1}^{3}}%
%BeginExpansion
{\displaystyle\sum\limits_{i=1}^{3}}
%EndExpansion
P_{p}^{\left(  0,0\right)  }\left(  1-2\lambda_{K,i}\right)  \right)  &
\text{on }K,\\
0 & \text{on }\Omega\backslash K.
\end{array}
\right.
\]
For $p$ odd and $E\in\mathcal{E}$, the \emph{non-conforming edge bubble} is
given by%
\[
B_{p,E}^{\operatorname*{CR}}:=\left\{
\begin{array}
[c]{ll}%
P_{p}^{\left(  0,0\right)  }\left(  1-2\lambda_{K,\mathbf{A}_{K,E}}\right)  &
\text{on }K\text{ for }K\in\mathcal{T}_{E},\\
0 & \text{on }\Omega\backslash\omega_{E},
\end{array}
\right.
\]
where $\mathbf{A}_{K,E}$ denotes the vertex in $K$ which is opposite to the
edge $E$.
\end{definition}

Different representations of the functions $B_{p,E}^{\operatorname*{CR}}$,
$B_{p,K}^{\operatorname*{CR}}$ exist in the literature, see \cite{BaranCVD},
\cite{Ainsworth_Rankin}, \cite[for $p=4,6.$]{ChaLeeLee}, \cite{ccss_2012}
while the formula for $B_{p,K}^{\operatorname*{CR}}$ has been introduced in
\cite{Baran_Stoyan}.

\begin{theorem}
\label{ThmBasisCRscalar}A basis for the space $CR_{p,0}\left(  \mathcal{T}%
\right)  $ is given

\begin{itemize}
\item for even $p$ by%
\[
\left\{  B_{p,\mathbf{z}},\mathbf{z}\in\mathcal{N}_{p,\Omega}\left(
\mathcal{T}\right)  \right\}  \cup\left\{  B_{p,K}^{\operatorname*{CR}}%
,K\in\mathcal{T}\right\}  ,
\]

\item for odd $p$ by%
\[
\left\{  B_{p,\mathbf{z}},\mathbf{z}\in\mathcal{N}_{p,\Omega}\left(
\mathcal{T}\right)  \backslash\mathcal{V}_{\Omega}\right\}  \cup\left\{
B_{p,E}^{\operatorname*{CR}},E\in\mathcal{E}_{\Omega}\right\}  .
\]

\end{itemize}
\end{theorem}

%

%TCIMACRO{\TeXButton{Proof}{\proof}}%
%BeginExpansion
\proof
%EndExpansion
For even $p$, this follows from \cite[Rem. 3]{BaranCVD} in combination with
\cite[Thm. 22]{ccss_2012}.

For odd $p$, we first observe that
\[
\left.  B_{p,E}^{\operatorname*{CR}}\right\vert _{E}=1\text{ so that }\left.
B_{p,E}^{\operatorname*{CR}}\right\vert _{\omega_{E}}\in C^{0}\left(
\omega_{E}\right)
\]
and for $K\in\mathcal{T}_{E}$ and $E^{\prime}\subset\partial K\backslash E$,
the restriction $\left.  B_{p,E}^{\operatorname*{CR}}\right\vert _{E^{\prime}%
}$ is the Legendre polynomial (lifted to the edge $E^{\prime}$) with endpoint
values $1$ at $\partial E$ and $-1$ at $\mathbf{A}_{K,E}$. Hence, the
assertion follows from \cite[Thm. 22]{ccss_2012}.%
%TCIMACRO{\TeXButton{End Proof}{\endproof}}%
%BeginExpansion
\endproof
%EndExpansion

\begin{corollary}
\label{CorBasis}A basis for the space $\mathbf{S}_{p,0}^{\operatorname*{CR}%
}\left(  \mathcal{T}\right)  $ is given

\begin{enumerate}
\item for even $p$ by%
\begin{equation}%
\begin{array}
[c]{l}%
\left\{  B_{p,\mathbf{z}}\mathbf{v}_{\mathbf{z}},\quad\mathbf{z}\in
\mathcal{N}_{p,\Omega}\left(  \mathcal{T}\right)  \right\}  \cup\left\{
B_{p,\mathbf{z}}\mathbf{w}_{\mathbf{z}},\quad\mathbf{z}\in\mathcal{N}%
_{p,\Omega}\left(  \mathcal{T}\right)  \right\} \\
\quad\\
\quad\cup\left\{  B_{p,K}^{\operatorname*{CR}}\mathbf{v}_{K},\quad
K\in\mathcal{T}\right\}  \cup\left\{  B_{p,K}^{\operatorname*{CR}}%
\mathbf{w}_{K},\quad K\in\mathcal{T}\right\}  ,
\end{array}
\label{basisvel}%
\end{equation}

\item for odd $p$ by%
\begin{equation}%
\begin{array}
[c]{l}%
\left\{  B_{p,\mathbf{z}}\mathbf{v}_{\mathbf{z}},\quad\mathbf{z}\in
\mathcal{N}_{p,\Omega}\left(  \mathcal{T}\right)  \backslash\mathcal{V}%
_{\Omega}\right\}  \cup\left\{  B_{p,\mathbf{z}}\mathbf{w}_{\mathbf{z}}%
,\quad\mathbf{z}\in\mathcal{N}_{p,\Omega}\left(  \mathcal{T}\right)
\backslash\mathcal{V}_{\Omega}\right\} \\
\quad\\
\quad\cup\left\{  B_{p,E}^{\operatorname*{CR}}\mathbf{v}_{E},\quad
E\in\mathcal{E}_{\Omega}\right\}  \cup\left\{  B_{p,E}^{\operatorname*{CR}%
}\mathbf{w}_{E},\quad E\in\mathcal{E}_{\Omega}\right\}  .
\end{array}
\label{basisvelodd}%
\end{equation}

\end{enumerate}

Here, for any nodal point $\mathbf{z}$, the linearly independent vectors
$\mathbf{v}_{\mathbf{z}},\mathbf{w}_{\mathbf{z}}\in\mathbb{R}^{2}$ can be
chosen arbitrarily. The same holds for any triangle $K$ for the vectors
$\mathbf{v}_{K},\mathbf{w}_{K}\in\mathbb{R}^{2}$ in (\ref{basisvel}) and for
any $E\in\mathcal{E}_{\Omega}$ for the vectors $\mathbf{v}_{E},\mathbf{w}%
_{E}\in\mathbb{R}^{2}$ in (\ref{basisvelodd}).
\end{corollary}

\begin{remark}
The original definition by \cite{CrouzeixRaviart} is implicit and given for
regular simplicial finite element meshes in $\mathbb{R}^{d}$, $d=2,3$. For
their practical implementation, a basis is needed and Corollary \ref{CorBasis}
provides a simple definition. A basis for Crouzeix-Raviart finite elements in
$\mathbb{R}^{3}$ has been introduced in \cite{Fortin_d3} for $p=2$ and a
general construction is given in \cite{CDS}.
\end{remark}

\subsection{The Pressure Kernel\label{PTpge5}}

For the investigation of the discrete inf-sup condition we employ the
macroelement technique in the form described in \cite{Stenberg_marco_new} (see
also \cite{Stenberg_macro}).

Let us first assume that every triangle $K\in\mathcal{T}$ has a vertex
$\mathbf{z}\in\Omega$. As a consequence, the sets $\mathcal{T}_{\mathbf{z}}$,
$\mathbf{z}\in\mathcal{V}_{\Omega}$, with nodal patches $\omega_{\mathbf{z}}$
form a macroelement partitioning of $\Omega$ in the sense of
\cite{Stenberg_marco_new}. We define the spaces%
\begin{align}
N_{p,\mathbf{z}}^{\operatorname*{CR}}  &  :=\left\{  p\in\mathbb{P}%
_{p-1}\left(  \mathcal{T}_{\mathbf{z}}\right)  \mid\forall\mathbf{v}%
\in\mathbf{CR}_{p,0}\left(  \mathcal{T}_{\mathbf{z}}\right)  :\left(
p,\operatorname*{div}\mathbf{v}\right)  _{L^{2}\left(  \omega_{\mathbf{z}%
}\right)  }=0\right\}  ,\label{defNCR}\\
N_{p,\mathbf{z}}  &  :=\left\{  p\in\mathbb{P}_{p-1}\left(  \mathcal{T}%
_{\mathbf{z}}\right)  \mid\forall\mathbf{v}\in\mathbf{S}_{p,0}\left(
\mathcal{T}_{\mathbf{z}}\right)  :\left(  p,\operatorname*{div}\mathbf{v}%
\right)  _{L^{2}\left(  \omega_{\mathbf{z}}\right)  }=0\right\}  .\nonumber
\end{align}

\begin{remark}
\label{Reminclusions}The definition of the Crouzeix-Raviart spaces
(\ref{PCR0RB}) implies $S_{p,0}\left(  \mathcal{T}_{\mathbf{z}}\right)
\subset CR_{p,0}\left(  \mathcal{T}_{\mathbf{z}}\right)  $ and, in turn,
\begin{equation}
N_{p,\mathbf{z}}^{\operatorname*{CR}}\subset N_{p,\mathbf{z}}. \label{inclN}%
\end{equation}
Let $1_{\mathbf{z}}:\omega_{\mathbf{z}}\rightarrow\mathbb{R}$ denote the
function with constant value $1$. Then, an integration by parts implies
$1_{\mathbf{z}}\in N_{p,\mathbf{z}}^{\operatorname*{CR}}$ so that $\dim
N_{p,\mathbf{z}}^{\operatorname*{CR}}\geq1$.
\end{remark}

The following Theorem is a direct consequence of \cite[Thm. 2.1]%
{Stenberg_marco_new}.

\begin{theorem}
\label{ThmStabQO}Let $\mathcal{T}$ be a regular finite element triangulation
of a bounded polygonal domain $\Omega\subset\mathbb{R}^{2}$ as in
\S \ref{NuDiscrete} with shape regularity constant $\gamma_{\mathcal{T}}$ and
at least one inner vertex. Let $p\in\mathbb{N}$. If
\begin{equation}
\dim N_{p,\mathbf{z}}^{\operatorname*{CR}}=1\qquad\forall\mathbf{z}%
\in\mathcal{V}_{\Omega}, \label{dimcond}%
\end{equation}
then the discrete inf-sup condition (\ref{infsupcond}) is satisfied with a
constant $c_{\mathcal{T},p}$ depending only on $p$, $\gamma_{\mathcal{T}}$,
and $d_{T}$ (cf. (\ref{defgammat}), (\ref{dtau})).

The discrete Stokes equation (\ref{discrStokes}) has a unique solution
$\left(  \mathbf{u}_{\operatorname*{CR}},p_{\operatorname*{disc}}\right)
\in\mathbf{S}_{p,0}^{\operatorname*{CR}}\left(  \mathcal{T}\right)  \times
M_{p-1}\left(  \mathcal{T}\right)  $.
\end{theorem}

\begin{remark}
If the assumptions of Theorem \ref{ThmStabQO} are satisfied, various types of
error estimates follow from well-established theory. Since this is not the
major theme of our paper, we briefly summarize two approaches: The error
$\left\Vert \mathbf{u}-\mathbf{u}_{\operatorname*{CR}}\right\Vert
_{\mathbf{H}^{1}\left(  \Omega\right)  }+\left\Vert p-p_{\operatorname*{disc}%
}\right\Vert _{L^{2}\left(  \Omega\right)  }$ can be estimated by means of the
second Strang lemma (see \cite{BergerScottStrang}, \cite[Thm. 4.2.2]%
{Ciarlet_orig}) via the sum of a quasi-optimal term and a consistency error.
The latter one converges with \emph{optimal rates} to zero by assuming
sufficiently high regularity of the continuous solution (see \cite[Thm.
3]{CrouzeixRaviart}, \cite[Thm. 2.2]{ChaLeeLee}).

There are also methods to establish \emph{quasi-optimality} of non-conforming
Crouzeix-Raviart discretizations. We mention the paper \cite{Veeser_qo}, where
a mapping from the non-conforming space to a conforming one is introduced and
employed in the discretization. For this method, quasi-optimal error estimates
can be proved.
\end{remark}

Assumption (\ref{dimcond}) is proved for $p=1,2$ and $p\geq4$, $p$ even. Also
note that we may restrict to triangulations with the property that any
$K\in\mathcal{T}$ has an inner vertex since the result for triangulations,
which contain at least one inner vertex, is implied by the following lemma.

\begin{lemma}
\label{Lemextension}Let $\mathcal{T}$ denote a regular triangulation of a
domain $\Omega$. Let $K^{\prime}$ denote a further triangle (not contained in
$\mathcal{T}$) which is attached to $\mathcal{T}$ by an edge $E$, more precisely:

\begin{enumerate}
\item $K^{\prime}\cap\Omega=\emptyset$ and there exists some $K\in\mathcal{T}$
and $E\in\mathcal{E}$ such that $K^{\prime}\cap K=E$,

\item Let $\mathcal{T}^{\prime}:=\mathcal{T\cup}\left\{  K\right\}  $ and
$\Omega^{\prime}:=\operatorname*{dom}\mathcal{T}^{\prime}$ (cf. (\ref{defdomT}%
)) is a bounded, polygonal Lipschitz domain in $\mathbb{R}^{2}$.
\end{enumerate}

If the pair $\mathbf{CR}_{p,0}\left(  \mathcal{T}\right)  \times
M_{p-1}\left(  \mathcal{T}\right)  $ is inf-sup stable, then $\mathbf{CR}%
_{p,0}\left(  \mathcal{T}^{\prime}\right)  \times M_{p-1}\left(
\mathcal{T}^{\prime}\right)  $ is inf-sup stable with a constant
$c_{\mathcal{T}^{\prime},p}$ which depends on $c_{\mathcal{T},p}$, the shape
regularity of $\left\{  K^{\prime}\right\}  $ and the polynomial degree $p$.
\end{lemma}

A proof for $p=3$ is given in \cite[Lem. 6.2]{Crouzeix_Falk}. Inspection of
the proof shows that it applies also to general $p$ provided the assumptions
of the lemma are satisfied. We omit the repetition of the arguments here.

Taking into account Lemma \ref{Lemextension} and the known cases from the
quoted literature, we assume for the following%
\begin{equation}%
\begin{array}
[c]{ll}%
\text{a)} & p\geq3\text{ is odd and}\\
\text{b)} & \mathcal{T}\text{ is a regular triangulation as in
\S \ref{NuDiscrete} s.t. every }K\in\mathcal{T}\text{ has one inner vertex.}%
\end{array}
\label{Assumptab}%
\end{equation}

In the following we investigate the condition (\ref{dimcond}) and start with
the definition of \textit{critical points }of a nodal patch $\mathcal{T}%
_{\mathbf{z}}$ for $\mathbf{z}\in\mathcal{V}_{\Omega}$ (see
\cite{ScottVogelius}).

\begin{definition}
\label{DefCritpoint}Let $\mathcal{T}$ denote a triangulation as in
\S \ref{NuDiscrete}. For $\mathbf{z}\in\mathcal{V}_{\Omega}$, we denote by
$\mathcal{T}_{\mathbf{z}}$, $\omega_{\mathbf{z}}$ the nodal patch as in
(\ref{nodalpatch}) and let
\[
\mathcal{E}_{\mathbf{z}}:=\left\{  E\in\mathcal{E}:E\subset\omega_{\mathbf{z}%
}\right\}  ,\qquad\mathcal{V}_{\mathbf{z}}:=\left\{  \mathbf{z}^{\prime}%
\in\mathcal{V}\mid\mathbf{z}^{\prime}\in\omega_{\mathbf{z}}\right\}  .
\]
A point $\mathbf{z}^{\prime}\in\mathcal{V}_{\mathbf{z}}$ is a \emph{critical
point} for $\mathcal{T}_{\mathbf{z}}$ if there exist two straight infinite
lines $L_{1}$, $L_{2}$ in $\mathbb{R}^{2}$ such that all edges $E\in
\mathcal{E}_{\mathbf{z}}$ having $\mathbf{z}^{\prime}$ as an endpoint satisfy
$E\subset L_{1}\cup L_{2}$. The set of all critical points in $\mathcal{T}%
_{\mathbf{z}}$ is $\mathcal{C}_{\mathbf{z}}$ and its cardinality denoted by
$\sigma_{\mathbf{z}}:=\left\vert \mathcal{C}_{\mathbf{z}}\right\vert $.
\end{definition}

\begin{figure}[ptb]%
\centering
\includegraphics[
height=2.1015in,
width=2.316in
]%
{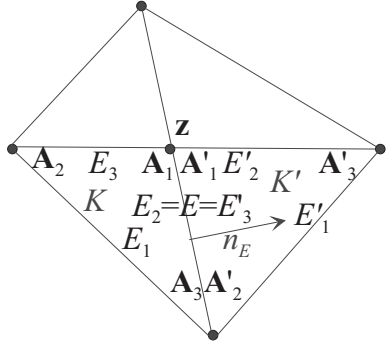}%
\caption{Geometric situation, where $\mathbf{z}$ is a critical vertex in the
patch $\mathcal{T}_{\mathbf{z}}$.}%
\label{Figcrisscross}%
\end{figure}

\begin{figure}[ptb]%
\centering
\includegraphics[
height=1.855in,
width=3.0355in
]%
{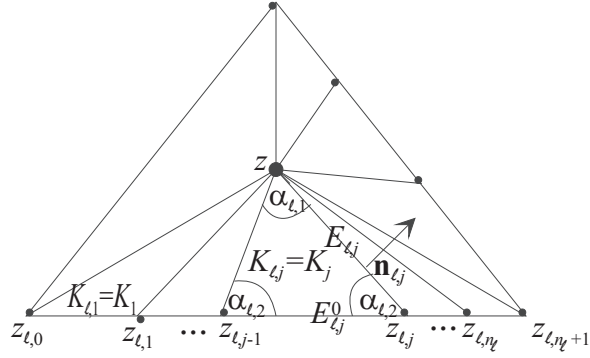}%
\caption{Nodal patch, illustrating singular vertices at $\partial
\omega_{\mathbf{z}}$. They are partitioned into edge-connected parts
$\mathcal{C}_{\mathbf{z},\ell}$, $1\leq\ell\leq c_{\mathbf{z}}$. The extremal
points $\mathbf{z}_{\ell,0}$ and $\mathbf{z}_{\ell,n_{\ell}+1}$ are no
singular points. The edge connecting $\mathbf{z}$ with $\mathbf{z}_{\ell,j}$
is denoted by $E_{\ell,j}$ and the connecting edge for the singular points
$\mathbf{z}_{\ell,j-1}$, $\mathbf{z}_{\ell,j}$ is denoted by $E_{\ell,j}^{0}$.
The normal vector $\mathbf{n}_{\ell,j}$ at $E_{\ell,j}$ is counterclockwise
oriented. Condition (\ref{finalcondcrpo}) is considered for all $\mathcal{C}%
_{\mathbf{z},\ell}$ separately and the index $\ell$ is skipped in the
computation.}%
\label{Fig_nod_patch}%
\end{figure}

\begin{remark}
\label{RemCritGeom}Geometric configurations where critical points occur are
well studied in the literature (see, e.g., \cite{ScottVogelius}). For nodal
patches $\mathcal{T}_{\mathbf{z}}$, any critical point belongs to one of the
following cases:

\begin{enumerate}
\item $\mathcal{T}_{\mathbf{z}}$ consists of four triangles and $\mathbf{z}$
is the intersections of the two diagonals in $\omega_{\mathbf{z}}$. Then
$\mathbf{z}$ is a critical point; see Fig. \ref{Figcrisscross}. In this case,
it holds $\sigma_{\mathbf{z}}=1$.

\item $\mathbf{z}^{\prime}\in\mathcal{V}_{\mathbf{z}}\backslash\left\{
\mathbf{z}\right\}  $. Let $E$ denote the edge with endpoints $\mathbf{z}$ and
$\mathbf{z}^{\prime}$ and let $K,K^{\prime}\in\mathcal{T}_{\mathbf{z}}$ be the
adjacent triangles. Then, $\mathbf{z}^{\prime}$ is a critical point if the sum
of the two angles at $\mathbf{z}^{\prime}$in the triangles $K,K^{\prime}%
\in\mathcal{T}_{\mathbf{z}}$ equals $\pi$; see Fig. \ref{Fig_nod_patch}. In
this case, it holds $\mathbf{z}\notin\mathcal{C}_{\mathbf{z}}$.
\end{enumerate}
\end{remark}

From \cite{ScottVogelius} we know that for $p\geq4$ it holds%
\begin{equation}
\dim N_{p,\mathbf{z}}=1+\sigma_{\mathbf{z}}. \label{dimform}%
\end{equation}
Remark \ref{Reminclusions} implies $1\leq\dim N_{p,\mathbf{z}}%
^{\operatorname*{CR}}\leq\dim N_{p,\mathbf{z}}$. In this section, we define
functions $q_{p-1,\mathbf{z}^{\prime}}$ for the critical points $\mathbf{z}%
^{\prime}\in\mathcal{C}_{\mathbf{z}}$ such that the $\left(  1+\sigma
_{\mathbf{z}}\right)  $ functions $1_{\mathbf{z}}$ and $q_{p-1,\mathbf{z}%
^{\prime}}$, $\mathbf{z}^{\prime}\in\mathcal{C}_{\mathbf{z}}$ are linearly
independent and belong to $N_{p,\mathbf{z}}$.

Afterwards, we will prove the implication:%
\begin{equation}
\left(  \eta=\beta_{0}1_{\mathbf{z}}+\sum_{\mathbf{z}^{\prime}\in
\mathcal{C}_{\mathbf{z}}}\beta_{\mathbf{z}^{\prime}}q_{p-1,\mathbf{z}^{\prime
}}\text{ satisfies }\eta\in N_{p,\mathbf{z}}^{\operatorname*{CR}}\right)
\implies\eta\in\operatorname*{span}\left\{  1_{\mathbf{z}}\right\}  .
\label{impl0}%
\end{equation}

Hence, in all cases in (\ref{Assumptab}), where the dimension formula
(\ref{dimform}) holds (e.g. for $p\geq4$) the condition (\ref{dimcond}) and,
in turn, the assumptions of Theorem \ref{ThmStabQO} are satisfied.\medskip

To prove (\ref{impl0}) it is sufficient to consider nodal patches with
critical points, i.e., $\mathcal{C}_{\mathbf{z}}\neq\emptyset$, and we will
construct basis functions for $N_{p,\mathbf{z}}$ explicitly. We fix a
(non-unique) sign function $\sigma:\mathcal{T}_{\mathbf{z}}\rightarrow\left\{
-1,1\right\}  $ by the condition:
\[
\text{if\ }K,K^{\prime}\in\mathcal{T}_{\mathbf{z}}\ \text{share an edge, then
}\sigma_{K}=-\sigma_{K^{\prime}}.
\]

\begin{definition}
\label{Defcritfunc}Let $\mathbf{z}^{\prime}\in\mathcal{C}_{\mathbf{z}}$ be a
critical point for $\mathcal{T}_{\mathbf{z}}$. The \emph{critical function}
$q_{p-1,\mathbf{z}^{\prime}}\in\mathbb{P}_{p-1}\left(  \mathcal{T}%
_{\mathbf{z}}\right)  $ for $\mathbf{z}^{\prime}$ is given by%
\[
q_{p-1,\mathbf{z}^{\prime}}:=\left\{
\begin{array}
[c]{ll}%
\frac{\sigma_{K}}{\left\vert K\right\vert }P_{p-1}^{\left(  0,2\right)
}\left(  1-2\lambda_{K,\mathbf{z}^{\prime}}\right)  & \text{on }%
K\in\mathcal{T}_{\mathbf{z}}\text{ with }\mathbf{z}^{\prime}\in K,\\
0 & \text{otherwise.}%
\end{array}
\right.
\]

\end{definition}

\begin{lemma}
\label{Lembasiscontker}Let (\ref{Assumptab}) be satisfied. The functions
$1_{\mathbf{z}}$ and $q_{p-1,\mathbf{z}^{\prime}}$, $\mathbf{z}^{\prime}%
\in\mathcal{C}_{\mathbf{z}}$, are linearly independent and belong to
$N_{p,\mathbf{z}}$. If the dimension formula (\ref{dimform}) holds, they form
a basis of $N_{p,\mathbf{z}}$.
\end{lemma}

%

%TCIMACRO{\TeXButton{Proof}{\proof}}%
%BeginExpansion
\proof
%EndExpansion
From Remark \ref{Reminclusions} it follows that $1_{z}\in N_{p,\mathbf{z}}$.
Next, we prove that $q_{p-1,\mathbf{z}^{\prime}}\in N_{p,\mathbf{z}}$. Let
$K\in\mathcal{T}_{\mathbf{z}}$ with $\mathbf{z}^{\prime}\in K$. Let
$\mathbf{A}_{i}$, $1\leq i\leq3$, denote the vertices of $K$ with the
convention $\mathbf{A}_{1}=\mathbf{z}$.

Recall the notation for barycentric coordinates in (\ref{lambdaintro1}),
(\ref{lambdaintro2}) and that $\mathbf{e}_{3,\ell}$, $\ell=1,2,3$, denotes the
$\ell$-th canonical unit vector in $\mathbb{R}^{3}$. Let $\mathbf{v}_{1}$,
$\mathbf{v}_{2}\in\mathbb{R}^{2}$ denote two linearly independent vectors --
the precise choice will be fixed later.

The restriction of the space $\mathbf{S}_{p,0}\left(  \mathcal{T}_{\mathbf{z}%
}\right)  $ to $K$ is spanned by (cf. (\ref{lambdaintro2}))%
\begin{equation}%
%TCIMACRO{\TeXButton{boldlambda}{\mbox{\boldmath$ \lambda$}}}%
%BeginExpansion
\mbox{\boldmath$ \lambda$}%
%EndExpansion
_{K}^{%
%TCIMACRO{\TeXButton{boldmu}{\mbox{\boldmath$ \mu$}}}%
%BeginExpansion
\mbox{\boldmath$ \mu$}%
%EndExpansion
}\mathbf{v}_{j}\qquad\text{for }j=1,2\text{ and }%
%TCIMACRO{\TeXButton{boldmu}{\mbox{\boldmath$ \mu$}}}%
%BeginExpansion
\mbox{\boldmath$ \mu$}%
%EndExpansion
=\mathbf{e}_{3,1}+%
%TCIMACRO{\TeXButton{boldalpha}{\mbox{\boldmath$ \alpha$}}}%
%BeginExpansion
\mbox{\boldmath$ \alpha$}%
%EndExpansion
\quad\text{with }%
%TCIMACRO{\TeXButton{boldalpha}{\mbox{\boldmath$ \alpha$}}}%
%BeginExpansion
\mbox{\boldmath$ \alpha$}%
%EndExpansion
\in\mathbb{I}_{=p-1}^{3}\text{.} \label{pointofvj}%
\end{equation}
Let $\mathbf{w}\in\left\{  \mathbf{v}_{1},\mathbf{v}_{2}\right\}  $ and let
$\partial_{\mathbf{w}}$ denote the G\^{a}teaux derivative of a function in the
direction $\mathbf{w}$. Then%
\[
\int_{K}q_{p-1,\mathbf{z}^{\prime}}\operatorname*{div}\left(
%TCIMACRO{\TeXButton{boldlambda}{\mbox{\boldmath$ \lambda$}}}%
%BeginExpansion
\mbox{\boldmath$ \lambda$}%
%EndExpansion
_{K}^{%
%TCIMACRO{\TeXButton{boldmu}{\mbox{\boldmath$ \mu$}}}%
%BeginExpansion
\mbox{\boldmath$ \mu$}%
%EndExpansion
}\mathbf{w}\right)  =\int_{K}q_{p-1,\mathbf{z}^{\prime}}\partial_{\mathbf{w}}%
%TCIMACRO{\TeXButton{boldlambda}{\mbox{\boldmath$ \lambda$}}}%
%BeginExpansion
\mbox{\boldmath$ \lambda$}%
%EndExpansion
_{K}^{%
%TCIMACRO{\TeXButton{boldmu}{\mbox{\boldmath$ \mu$}}}%
%BeginExpansion
\mbox{\boldmath$ \mu$}%
%EndExpansion
}=\sum_{s=1}^{3}\mu_{s}\partial_{\mathbf{w}}\lambda_{K,s}\int_{K}%
q_{p-1,\mathbf{z}^{\prime}}%
%TCIMACRO{\TeXButton{boldlambda}{\mbox{\boldmath$ \lambda$}}}%
%BeginExpansion
\mbox{\boldmath$ \lambda$}%
%EndExpansion
_{K}^{%
%TCIMACRO{\TeXButton{boldmu}{\mbox{\boldmath$ \mu$}}}%
%BeginExpansion
\mbox{\boldmath$ \mu$}%
%EndExpansion
-\mathbf{e}_{3,s}}.
\]
There exists $\ell\in\left\{  1,2,3\right\}  $ such that $\mathbf{z}^{\prime
}=\mathbf{A}_{\ell}$. We employ \cite[Prop. 4.1]{Baran_Stoyan} which tells us
(as a consequence of \cite[(3.10)]{Koornwinder_75}) that for $j\in\left\{
1,2,3\right\}  \backslash\left\{  \ell\right\}  $ it holds%
\begin{equation}
\int_{K}w\lambda_{j}P_{p-1}^{\left(  0,2\right)  }\left(  1-2\lambda_{K,\ell
}\right)  =0\quad\forall w\in\mathbb{P}_{p-2}\left(  K\right)  .
\label{orthoconc}%
\end{equation}
For $j\in\left\{  1,2,3\right\}  \backslash\left\{  \ell\right\}  $ with
$\mu_{j}-\delta_{j,s}\geq1$, this implies%
\begin{equation}
I_{K}^{\left(  \ell,s\right)  }:=\mu_{s}\int_{K}q_{p-1,\mathbf{A}_{\ell}}%
%TCIMACRO{\TeXButton{boldlambda}{\mbox{\boldmath$ \lambda$}}}%
%BeginExpansion
\mbox{\boldmath$ \lambda$}%
%EndExpansion
_{K}^{%
%TCIMACRO{\TeXButton{boldmu}{\mbox{\boldmath$ \mu$}}}%
%BeginExpansion
\mbox{\boldmath$ \mu$}%
%EndExpansion
-\mathbf{e}_{3,s}}=0. \label{muesI}%
\end{equation}
Taking into account the orthogonality relations (\ref{orthoconc}) and that the
factor $\mu_{s}$ in (\ref{muesI}) vanishes in many cases, it remains to
consider the cases%
\begin{align}
\left(  \ell,s\right)   &  =\left(  1,1\right)  \text{ and }%
%TCIMACRO{\TeXButton{boldmu}{\mbox{\boldmath$ \mu$}}}%
%BeginExpansion
\mbox{\boldmath$ \mu$}%
%EndExpansion
=\left(  p,0,0\right)  ,\nonumber\\
\left(  \ell,s\right)   &  =\left(  2,1\right)  \text{ and }%
%TCIMACRO{\TeXButton{boldmu}{\mbox{\boldmath$ \mu$}}}%
%BeginExpansion
\mbox{\boldmath$ \mu$}%
%EndExpansion
=\left(  1,p-1,0\right)  ,\nonumber\\
\left(  \ell,s\right)   &  =\left(  3,1\right)  \text{ and }%
%TCIMACRO{\TeXButton{boldmu}{\mbox{\boldmath$ \mu$}}}%
%BeginExpansion
\mbox{\boldmath$ \mu$}%
%EndExpansion
=\left(  1,0,p-1\right)  ,\label{lsremcase}\\
\left(  \ell,s\right)   &  =\left(  1,2\right)  \text{ and }%
%TCIMACRO{\TeXButton{boldmu}{\mbox{\boldmath$ \mu$}}}%
%BeginExpansion
\mbox{\boldmath$ \mu$}%
%EndExpansion
=\left(  p-1,1,0\right)  ,\nonumber\\
\left(  \ell,s\right)   &  =\left(  1,3\right)  \text{ and }%
%TCIMACRO{\TeXButton{boldmu}{\mbox{\boldmath$ \mu$}}}%
%BeginExpansion
\mbox{\boldmath$ \mu$}%
%EndExpansion
=\left(  p-1,0,1\right) \nonumber
\end{align}
with corresponding integrals%
\[%
\begin{array}
[c]{lll}%
I_{K}^{\left(  1,1\right)  }=p\int_{K}q_{p-1,\mathbf{A}_{1}}\lambda
_{K,1}^{p-1}, & I_{K}^{\left(  2,1\right)  }=\int_{K}q_{p-1,\mathbf{A}_{2}%
}\lambda_{K,2}^{p-1}, & I_{K}^{\left(  3,1\right)  }=\int_{K}q_{p-1,\mathbf{A}%
_{3}}\lambda_{K,3}^{p-1},\\
I_{K}^{\left(  1,2\right)  }=\int_{K}q_{p-1,\mathbf{A}_{1}}\lambda_{K,1}%
^{p-1}, & I_{K}^{\left(  1,3\right)  }=\int_{K}q_{p-1,\mathbf{A}_{1}}%
\lambda_{K,1}^{p-1}. &
\end{array}
\]
All other integrals in (\ref{muesI}) vanish.

A standard affine transformation to the reference element $\widehat{K}$ shows
that%
\begin{equation}
\frac{1}{p}I_{K}^{\left(  1,1\right)  }=I_{K}^{\left(  2,1\right)  }%
=I_{K}^{\left(  3,1\right)  }=I_{K}^{\left(  1,2\right)  }=I_{K}^{\left(
1,3\right)  }=\sigma_{K}c_{p}, \label{Ipequal}%
\end{equation}
where%
\[
c_{p}:=2\int_{\widehat{K}}P_{p-1}^{\left(  0,2\right)  }\left(
1-2\widehat{\lambda}_{2}\right)  \widehat{\lambda}_{2}^{p-1}.
\]
It remains to prove the second equality in%
\[
\int_{K}q_{p-1,\mathbf{z}^{\prime}}\operatorname*{div}\left(
%TCIMACRO{\TeXButton{boldlambda}{\mbox{\boldmath$ \lambda$}}}%
%BeginExpansion
\mbox{\boldmath$ \lambda$}%
%EndExpansion
_{K}^{%
%TCIMACRO{\TeXButton{boldmu}{\mbox{\boldmath$ \mu$}}}%
%BeginExpansion
\mbox{\boldmath$ \mu$}%
%EndExpansion
}\mathbf{w}\right)  =\sum_{s=1}^{3}I_{K}^{\left(  \ell,s\right)  }%
\partial_{\mathbf{w}}\lambda_{K,s}=0
\]
for $\left(  \ell,s\right)  \in\left\{  \left(  1,r\right)  ,\left(
r,1\right)  :r\in\left\{  1,2,3\right\}  \right\}  $ (cf. (\ref{lsremcase})).
By using (\ref{Ipequal}) we obtain%
\begin{equation}
\int_{K}q_{p-1,\mathbf{A}_{\ell}}\operatorname*{div}\left(
%TCIMACRO{\TeXButton{boldlambda}{\mbox{\boldmath$ \lambda$}}}%
%BeginExpansion
\mbox{\boldmath$ \lambda$}%
%EndExpansion
_{K}^{%
%TCIMACRO{\TeXButton{boldmu}{\mbox{\boldmath$ \mu$}}}%
%BeginExpansion
\mbox{\boldmath$ \mu$}%
%EndExpansion
}\mathbf{w}\right)  =\sigma_{K}c_{p}\left\{
\begin{array}
[c]{ll}%
p\partial_{\mathbf{w}}\lambda_{K,1} & \ell=1\wedge%
%TCIMACRO{\TeXButton{boldmu}{\mbox{\boldmath$ \mu$}}}%
%BeginExpansion
\mbox{\boldmath$ \mu$}%
%EndExpansion
=\left(  p,0,0\right)  ,\\
\partial_{\mathbf{w}}\lambda_{K,1} & \ell=2\wedge%
%TCIMACRO{\TeXButton{boldmu}{\mbox{\boldmath$ \mu$}}}%
%BeginExpansion
\mbox{\boldmath$ \mu$}%
%EndExpansion
=\left(  1,p-1,0\right)  ,\\
\partial_{\mathbf{w}}\lambda_{K,1} & \ell=3\wedge%
%TCIMACRO{\TeXButton{boldmu}{\mbox{\boldmath$ \mu$}}}%
%BeginExpansion
\mbox{\boldmath$ \mu$}%
%EndExpansion
=\left(  1,0,p-1\right)  ,\\
\partial_{\mathbf{w}}\lambda_{K,2} & \ell=1\wedge%
%TCIMACRO{\TeXButton{boldmu}{\mbox{\boldmath$ \mu$}}}%
%BeginExpansion
\mbox{\boldmath$ \mu$}%
%EndExpansion
=\left(  p-1,1,0\right)  ,\\
\partial_{\mathbf{w}}\lambda_{K,3} & \ell=1\wedge%
%TCIMACRO{\TeXButton{boldmu}{\mbox{\boldmath$ \mu$}}}%
%BeginExpansion
\mbox{\boldmath$ \mu$}%
%EndExpansion
=\left(  p-1,0,1\right)  ,\\
0 & \text{otherwise.}%
\end{array}
\right.  \label{contelmat}%
\end{equation}

\textbf{Evaluation of the right-hand side in (\ref{contelmat}).}

The basis functions in $\mathbf{S}_{p,0}\left(  \mathcal{T}_{\mathbf{z}%
}\right)  $ can be grouped into three different types.

\begin{enumerate}
\item Basis functions whose support is one triangle $K\in\mathcal{T}%
_{\mathbf{z}}$. Let $\mathbf{A}_{i}$, $1\leq i\leq3$, denote the vertices of
$K$. Then, these basis functions are given by%
\[
B_{K}^{%
%TCIMACRO{\TeXButton{boldmu}{\mbox{\boldmath$ \mu$}}}%
%BeginExpansion
\mbox{\boldmath$ \mu$}%
%EndExpansion
,j}:=%
%TCIMACRO{\TeXButton{boldlambda}{\mbox{\boldmath$ \lambda$}}}%
%BeginExpansion
\mbox{\boldmath$ \lambda$}%
%EndExpansion
_{K}^{%
%TCIMACRO{\TeXButton{boldmu}{\mbox{\boldmath$ \mu$}}}%
%BeginExpansion
\mbox{\boldmath$ \mu$}%
%EndExpansion
}\mathbf{v}_{j},\qquad\text{for }%
%TCIMACRO{\TeXButton{boldmu}{\mbox{\boldmath$ \mu$}}}%
%BeginExpansion
\mbox{\boldmath$ \mu$}%
%EndExpansion
\in\mathbb{I}_{=p}^{3}\text{ and }\min\left\{  \mu_{i}:1\leq i\leq3\right\}
\geq1.
\]
From (\ref{contelmat}) it follows that%
\[
\int_{\omega_{\mathbf{z}}}q_{p-1,\mathbf{z}^{\prime}}\operatorname*{div}%
\left(  B_{K}^{%
%TCIMACRO{\TeXButton{boldmu}{\mbox{\boldmath$ \mu$}}}%
%BeginExpansion
\mbox{\boldmath$ \mu$}%
%EndExpansion
,j}\right)  =0.
\]

\item Basis functions whose support is an edge patch. Let $E\in\mathcal{E}%
_{\mathbf{z}}$ with $\mathbf{z\in}E$ and denote by $K$ and $K^{\prime}$ the
triangles in $\mathcal{T}_{\mathbf{z}}$ which share $E$. We denote the
vertices in $K$ and $K^{\prime}$ by $\mathbf{A}_{\ell}\ $and $\mathbf{A}%
_{\ell}^{\prime}$, $\ell=1,2,3$, with the convention that $\mathbf{A}%
_{1}=\mathbf{A}_{1}^{\prime}=\mathbf{z}$, $\mathbf{A}_{2}$ is the vertex in
$K$ opposite to $E$, and $\mathbf{A}_{3}^{\prime}$ is the vertex in
$K^{\prime}$ opposite to $E$. Thus, $\mathbf{A}_{3}=\mathbf{A}_{2}^{\prime}$.
Then, $2\left(  p-1\right)  $ basis functions of $\mathbf{S}_{p,0}\left(
\mathcal{T}_{\mathbf{z}}\right)  $ are associated with the edge $E$. For $%
%TCIMACRO{\TeXButton{boldmu}{\mbox{\boldmath$ \mu$}}}%
%BeginExpansion
\mbox{\boldmath$ \mu$}%
%EndExpansion
\in\mathbb{I}_{=p}^{3}$ with $%
%TCIMACRO{\TeXButton{boldmu}{\mbox{\boldmath$ \mu$}}}%
%BeginExpansion
\mbox{\boldmath$ \mu$}%
%EndExpansion
^{\prime}=\left(  \mu_{1},\mu_{3},\mu_{2}\right)  $, define%
\begin{equation}
B_{E}^{%
%TCIMACRO{\TeXButton{boldmu}{\mbox{\boldmath$ \mu$}}}%
%BeginExpansion
\mbox{\boldmath$ \mu$}%
%EndExpansion
,j}:=\left\{
\begin{array}
[c]{ll}%
%TCIMACRO{\TeXButton{boldlambda}{\mbox{\boldmath$ \lambda$}}}%
%BeginExpansion
\mbox{\boldmath$ \lambda$}%
%EndExpansion
_{K}^{%
%TCIMACRO{\TeXButton{boldmu}{\mbox{\boldmath$ \mu$}}}%
%BeginExpansion
\mbox{\boldmath$ \mu$}%
%EndExpansion
}\mathbf{v}_{j} & \text{on }K,\\%
%TCIMACRO{\TeXButton{boldlambda}{\mbox{\boldmath$ \lambda$}}}%
%BeginExpansion
\mbox{\boldmath$ \lambda$}%
%EndExpansion
_{K^{\prime}}^{%
%TCIMACRO{\TeXButton{boldmu}{\mbox{\boldmath$ \mu$}}}%
%BeginExpansion
\mbox{\boldmath$ \mu$}%
%EndExpansion
^{\prime}}\mathbf{v}_{j} & \text{on }K^{\prime},\\
0 & \text{on }\omega_{\mathbf{z}}\backslash\omega_{E}%
\end{array}
\right.  \quad\text{for }\left\{
\begin{array}
[c]{ll}
&
%TCIMACRO{\TeXButton{boldmu}{\mbox{\boldmath$ \mu$}}}%
%BeginExpansion
\mbox{\boldmath$ \mu$}%
%EndExpansion
\in\mathbb{I}_{=p}^{3},\\
\wedge & \min\left\{  \mu_{1},\mu_{3}\right\}  \geq1,\\
\wedge & \mu_{2}=0.
\end{array}
\right.  \label{oneedge}%
\end{equation}
The integral $\int_{\omega_{\mathbf{z}}}q_{p-1,\mathbf{z}^{\prime}%
}\operatorname*{div}\left(  B_{E}^{%
%TCIMACRO{\TeXButton{boldmu}{\mbox{\boldmath$ \mu$}}}%
%BeginExpansion
\mbox{\boldmath$ \mu$}%
%EndExpansion
,j}\right)  $ can be different from zero, only if $\mathbf{z}^{\prime}%
\in\left\{  \mathbf{A}_{1},\mathbf{A}_{2},\mathbf{A}_{3},\mathbf{A}%
_{3}^{\prime}\right\}  $.

\begin{enumerate}
\item If $\mathbf{z}^{\prime}=\mathbf{A}_{2}$, we get%
\[
\int_{\omega_{\mathbf{z}}}q_{p-1,\mathbf{z}^{\prime}}\operatorname*{div}%
\left(  B_{E}^{%
%TCIMACRO{\TeXButton{boldmu}{\mbox{\boldmath$ \mu$}}}%
%BeginExpansion
\mbox{\boldmath$ \mu$}%
%EndExpansion
,j}\right)  =\int_{K}q_{p-1,\mathbf{A}_{2}}\operatorname*{div}\left(
%TCIMACRO{\TeXButton{boldlambda}{\mbox{\boldmath$ \lambda$}}}%
%BeginExpansion
\mbox{\boldmath$ \lambda$}%
%EndExpansion
_{K}^{%
%TCIMACRO{\TeXButton{boldmu}{\mbox{\boldmath$ \mu$}}}%
%BeginExpansion
\mbox{\boldmath$ \mu$}%
%EndExpansion
}\mathbf{v}_{j}\right)  .
\]
The combination of (\ref{contelmat}) with the conditions for $%
%TCIMACRO{\TeXButton{boldmu}{\mbox{\boldmath$ \mu$}}}%
%BeginExpansion
\mbox{\boldmath$ \mu$}%
%EndExpansion
$ in (\ref{oneedge}) and $\ell=2$ (since $\mathbf{z}^{\prime}=\mathbf{A}_{2}$)
shows that these integrals vanish. We argue along the same lines to see that
the integrals vanish if $\mathbf{z}^{\prime}=\mathbf{A}_{3}^{\prime}$.

\item Let $\mathbf{z}^{\prime}=\mathbf{A}_{1}$. Then,%
\[
\int_{\omega_{E}}q_{p-1,\mathbf{z}^{\prime}}\operatorname*{div}\left(  B_{E}^{%
%TCIMACRO{\TeXButton{boldmu}{\mbox{\boldmath$ \mu$}}}%
%BeginExpansion
\mbox{\boldmath$ \mu$}%
%EndExpansion
,j}\right)  =\int_{K}q_{p-1,\mathbf{A}_{1}}\operatorname*{div}\left(
%TCIMACRO{\TeXButton{boldlambda}{\mbox{\boldmath$ \lambda$}}}%
%BeginExpansion
\mbox{\boldmath$ \lambda$}%
%EndExpansion
_{K}^{%
%TCIMACRO{\TeXButton{boldmu}{\mbox{\boldmath$ \mu$}}}%
%BeginExpansion
\mbox{\boldmath$ \mu$}%
%EndExpansion
}\mathbf{v}_{j}\right)  +\int_{K^{\prime}}q_{p-1,\mathbf{A}_{1}}%
\operatorname*{div}\left(
%TCIMACRO{\TeXButton{boldlambda}{\mbox{\boldmath$ \lambda$}}}%
%BeginExpansion
\mbox{\boldmath$ \lambda$}%
%EndExpansion
_{K^{\prime}}^{%
%TCIMACRO{\TeXButton{boldmu}{\mbox{\boldmath$ \mu$}}}%
%BeginExpansion
\mbox{\boldmath$ \mu$}%
%EndExpansion
^{\prime}}\mathbf{v}_{j}\right)  .
\]
From (\ref{contelmat}) with the conditions for $%
%TCIMACRO{\TeXButton{boldmu}{\mbox{\boldmath$ \mu$}}}%
%BeginExpansion
\mbox{\boldmath$ \mu$}%
%EndExpansion
$ in (\ref{oneedge}) and $\ell=1$ (since $\mathbf{z}^{\prime}=\mathbf{A}_{1}%
$), we see that these integrals vanish if $%
%TCIMACRO{\TeXButton{boldmu}{\mbox{\boldmath$ \mu$}}}%
%BeginExpansion
\mbox{\boldmath$ \mu$}%
%EndExpansion
\neq\left(  p-1,0,1\right)  $. However, for the case $%
%TCIMACRO{\TeXButton{boldmu}{\mbox{\boldmath$ \mu$}}}%
%BeginExpansion
\mbox{\boldmath$ \mu$}%
%EndExpansion
=\left(  p-1,0,1\right)  $ we conclude from (\ref{contelmat}) that%
\begin{equation}
\int_{\omega_{\mathbf{z}}}q_{p-1,\mathbf{z}^{\prime}}\operatorname*{div}%
\left(  B_{E}^{%
%TCIMACRO{\TeXButton{boldmu}{\mbox{\boldmath$ \mu$}}}%
%BeginExpansion
\mbox{\boldmath$ \mu$}%
%EndExpansion
,j}\right)  =\sigma_{K}c_{p}\left(  \partial_{\mathbf{v}_{j}}\lambda
_{K,\mathbf{A}_{3}}-\partial_{\mathbf{v}_{j}}\lambda_{K^{\prime}%
,\mathbf{A}_{2}^{\prime}}\right)  . \label{diffedfgecase}%
\end{equation}
Since $\mathbf{z}^{\prime}=\mathbf{A}_{1}=\mathbf{A}_{1}^{\prime}$ is a
critical point, the sum of the angles in $K$ at $\mathbf{z}^{\prime}$ and in
$K^{\prime}$ at $\mathbf{z}^{\prime}$ is $\pi$. Hence, $K\cup K^{\prime
}=\omega_{E}$ is the triangle with vertices $\mathbf{A}_{2},\mathbf{A}_{3}$,
$\mathbf{A}_{3}^{\prime}$ and $\left.  \lambda_{\omega_{E},\mathbf{A}_{3}%
}\right\vert _{K}=\lambda_{K,\mathbf{A}_{3}}$ and $\left.  \lambda_{\omega
_{E},\mathbf{A}_{3}}\right\vert _{K^{\prime}}=\lambda_{K^{\prime}%
,\mathbf{A}_{2}^{\prime}}$. Since $\lambda_{\omega_{E},\mathbf{A}_{3}}$ is
affine on $\omega_{E}$, the difference in (\ref{diffedfgecase}) vanishes.
\end{enumerate}

Let $\mathbf{z}^{\prime}=\mathbf{A}_{3}$. By arguing as in the previous case,
we obtain $\int_{\omega_{\mathbf{z}}}q_{p-1,\mathbf{z}^{\prime}}%
\operatorname*{div}\left(  B_{E}^{%
%TCIMACRO{\TeXButton{boldmu}{\mbox{\boldmath$ \mu$}}}%
%BeginExpansion
\mbox{\boldmath$ \mu$}%
%EndExpansion
,j}\right)  =0$.

\item Finally, it remains to consider the basis functions associated with the
centre node. We set%
\[
\left.  B_{\mathbf{z}}^{p,j}\right\vert _{K}:=\lambda_{K,\mathbf{A}_{1}}%
^{p}\mathbf{v}_{j}\quad\text{for }K\in\mathcal{T}_{\mathbf{z}}%
\]
with the numbering convention that the first vertex in $K$ satisfies
$\mathbf{z}=\mathbf{A}_{1}$.

\textbf{Case a. }Let the critical point be an \textquotedblleft
outer\textquotedblright\ vertex, i.e., $\mathbf{z}^{\prime}\in\mathcal{V}%
_{\mathbf{z}}\backslash\left\{  \mathbf{z}\right\}  $. The edge connecting
$\mathbf{z}$ with $\mathbf{z}^{\prime}$ is denoted by $E$ with adjacent
triangles $K$, $K^{\prime}$ and we employ the numbering convening as in case
2. Then,%
\[
\int_{\omega_{\mathbf{z}}}q_{p-1,\mathbf{z}^{\prime}}\operatorname*{div}%
\left(  B_{\mathbf{z}}^{p,j}\right)  =\int_{K}q_{p-1,\mathbf{A}_{3}%
}\operatorname*{div}\left(  \lambda_{K,\mathbf{A}_{1}}^{p}\mathbf{v}%
_{j}\right)  +\int_{K^{\prime}}q_{p-1,\mathbf{A}_{2}^{\prime}}%
\operatorname*{div}\left(  \lambda_{K^{\prime},\mathbf{A}_{1}^{\prime}}%
^{p}\mathbf{v}_{j}\right)  \overset{\text{(\ref{contelmat})}}{=}0.
\]

\textbf{Case b. }Let $\mathbf{z}^{\prime}=\mathbf{z}$ be a critical point.
From Remark \ref{RemCritGeom} we know that in this case the nodal patch
consists of $4$ triangles so that $\omega_{\mathbf{z}}$ is a quadrilateral and
$\mathbf{z}$ is the crosspoint of the diagonals $d_{1}$, $d_{2}$ in
$\omega_{\mathbf{z}}$. Fix $\mathbf{v}_{1}$ and $\mathbf{v}_{2}$ in
(\ref{pointofvj}) as a unit vector pointing in the direction of $d_{1}$ and
$d_{2}$.

We first consider the basis function $B_{\mathbf{z}}^{p,1}\ $corresponding to
the vector $\mathbf{v}_{1}$ aligned with $d_{1}$. We choose an edge
$E\in\mathcal{E}_{\mathbf{z}}$ with $\mathbf{z}\in E$ and $E\subset d_{1}$
such that $\mathbf{v}_{1}$ is tangential to $E$. The two triangles in
$\mathcal{T}_{\mathbf{z}}$ adjacent to $E$ are denoted by $K,K^{\prime}$ with
the same numbering convention of the vertices as in case 2. Then,%
\[
\int_{K\cup K^{\prime}}q_{p-1,\mathbf{z}}\operatorname*{div}\left(
B_{\mathbf{z}}^{p,1}\right)  =\int_{K}q_{p-1,\mathbf{A}_{1}}%
\operatorname*{div}\left(  \lambda_{K,\mathbf{A}_{1}}^{p}\mathbf{v}%
_{1}\right)  +\int_{K^{\prime}}q_{p-1,\mathbf{A}_{1}}\operatorname*{div}%
\left(  \lambda_{K^{\prime},\mathbf{A}_{1}}^{p}\mathbf{v}_{1}\right)  .
\]
From (\ref{contelmat}) we conclude that%
\begin{equation}
\int_{K\cup K^{\prime}}q_{p-1,\mathbf{z}}\operatorname*{div}\left(
B_{\mathbf{z}}^{p,1}\right)  =pc_{p}\sigma_{K}\left(  \partial_{\mathbf{v}%
_{1}}\lambda_{K,\mathbf{z}}-\partial_{\mathbf{v}_{1}}\lambda_{K^{\prime
},\mathbf{z}}\right)  . \label{baryfinal}%
\end{equation}
Since $\mathbf{v}_{1}$ is tangential to $E$ and the jump of $\lambda
_{K,\mathbf{z}}$ and $\lambda_{K^{\prime},\mathbf{z}}$ across $E$ is zero, the
same holds for the tangential derivative and the difference in
(\ref{baryfinal}) is zero.

For the other edge in $\mathcal{E}_{\mathbf{z}}$, which contains $\mathbf{z}$
and is aligned to $\mathbf{v}_{1}$, we argue in the same way to see that the
integral over the remaining two triangles in $\mathcal{T}_{\mathbf{z}}$ is
also zero,%
\[
\int_{\omega_{\mathbf{z}}}q_{p-1,\mathbf{z}}\operatorname*{div}\left(
B_{\mathbf{z}}^{p,1}\right)  =0.
\]

For the basis function $B_{\mathbf{z}}^{p,2}$ we consider the two edges
$E\in\mathcal{E}_{\mathbf{z}}$ with $\mathbf{z}\in E$ and which are aligned to
$\mathbf{v}_{2}$. By repeating the previous arguments we conclude
$\int_{\omega_{\mathbf{z}}}q_{p-1,\mathbf{z}}\operatorname*{div}\left(
B_{\mathbf{z}}^{p,2}\right)  =0$.
\end{enumerate}

In summary, we have proved that%
\[
1_{\mathbf{z}}\in N_{p,\mathbf{z}}\quad\text{and\quad}q_{p-1,\mathbf{z}%
^{\prime}}\in N_{p,\mathbf{z}}%
\]
for all critical points $\mathbf{z}^{\prime}\in\mathcal{C}_{\mathbf{z}}%
$.\medskip

It remains to prove that these functions are linearly independent.

If $\mathbf{z}^{\prime}=\mathbf{z}$ is a critical point in the interior of a
patch, we have the geometric situation as described in case 3b and no further
critical point exists. Since for $p\geq2$, $P_{p-1}^{\left(  0,2\right)
}\left(  1-2\lambda_{K,\mathbf{z}^{\prime}}\right)  $ is not constant, it
follows that $1_{\mathbf{z}}$ and $q_{p-1,\mathbf{z}^{\prime}}$ are linearly
independent.\medskip

It remains to consider the case that all critical points $\mathbf{z}^{\prime}$
in $\mathcal{T}_{z}$ are located on $\partial\omega_{\mathbf{z}}$. Then, the
support of the corresponding critical function $q_{p-1,\mathbf{z}^{\prime}}$
consists of only those two triangles which share the edge with endpoints
$\mathbf{z}$ and $\mathbf{z}^{\prime}$. Hence, it is sufficient to prove that
for a triangle $K\in\mathcal{T}_{\mathbf{z}}$ with vertices $\mathbf{A}_{1}$,
$\mathbf{A}_{2}$, $\mathbf{A}_{3}$ (convention: $\mathbf{z}=\mathbf{A}_{1}$)
and the property that both vertices $\mathbf{A}_{2}$ and $\mathbf{A}_{3}$ are
critical points it holds
\[
P_{p-1}^{\left(  0,2\right)  }\left(  1-2\lambda_{K,\mathbf{A}_{2}}\right)
\text{,\quad}P_{p-1}^{\left(  0,2\right)  }\left(  1-2\lambda_{K,\mathbf{A}%
_{3}}\right)  ,\quad1_{K}%
\]
are linearly independent, where $1_{K}$ is the constant function $1$ on $K$.
We transfer this claim to the reference element by an affine transformation
and see that it suffices to prove that%
\[
\hat{q}_{1}:=P_{p-1}^{\left(  0,2\right)  }\left(  1-2\widehat{\lambda}%
_{2}\right)  \text{,\quad}\hat{q}_{2}:=P_{p-1}^{\left(  0,2\right)  }\left(
1-2\widehat{\lambda}_{3}\right)  ,\quad\hat{q}_{3}:=1_{\widehat{K}}%
\]
are linearly independent. We will show that Gram's matrix $\mathbf{G}=\left(
g_{i,j}\right)  _{i,j=1}^{3}$ with $g_{i,j}:=\left(  \hat{q}_{i},\hat{q}%
_{j}\right)  _{L^{2}\left(  \widehat{K}\right)  }$ is regular. We use%
\[
P_{p-1}^{\left(  0,2\right)  }\left(  1-2\widehat{\lambda}_{3}\right)
=P_{p-1}^{\left(  0,2\right)  }\left(  1\right)  +\widehat{\lambda}_{3}%
r_{p-2}\quad\text{with\quad}r_{p-2}\in\mathbb{P}_{p-2}\left(  \left[
-1,1\right]  \right)  \text{.}%
\]
Hence, the combination of (\ref{orthoconc}) with $P_{p-1}^{\left(  0,2\right)
}\left(  1\right)  =1$ (cf. \cite[Table 18.6.1]{NIST:DLMF}) leads to%
\[
g_{12}=g_{21}=d_{p}P_{p-1}^{\left(  0,2\right)  }\left(  1\right)  =d_{p},
\]
where
\begin{align*}
d_{p}  &  :=\int_{\hat{K}}P_{p-1}^{\left(  0,2\right)  }\left(
1-2\widehat{\lambda}_{2}\right)  =\int_{\hat{K}}P_{p-1}^{\left(  0,2\right)
}\left(  1-2\widehat{\lambda}_{3}\right) \\
&  =\frac{1}{4}\int_{-1}^{1}\left(  1+t\right)  P_{p-1}^{\left(  0,2\right)
}\left(  t\right)  dt\overset{\text{Lem. \ref{LemCHJP}}}{=}\frac{\left(
-1\right)  ^{p-1}}{p\left(  p+1\right)  }.
\end{align*}
Furthermore, it holds%
\begin{align*}
g_{11}  &  =g_{22}=\int_{\hat{K}}\left(  P_{p-1}^{\left(  0,2\right)  }\left(
1-2\widehat{\lambda}_{2}\right)  \right)  ^{2}=\int_{0}^{1}\left(
P_{p-1}^{\left(  0,2\right)  }\left(  1-2x_{1}\right)  \right)  ^{2}\left(
1-x_{1}\right)  dx_{1}\\
&  =\frac{1}{4}\int_{-1}^{1}\left(  P_{p-1}^{\left(  0,2\right)  }\left(
t\right)  \right)  ^{2}\left(  1+t\right)  dt\overset{\text{Lem.
\ref{LemCHJP}}}{=}\frac{1}{2}.
\end{align*}
Since $p$ is odd we deduce%
\[
\mathbf{G}=\left[
\begin{array}
[c]{lll}%
\frac{1}{2} & \frac{1}{p\left(  p+1\right)  } & \frac{1}{p\left(  p+1\right)
}\\
\frac{1}{p\left(  p+1\right)  } & \frac{1}{2} & \frac{1}{p\left(  p+1\right)
}\\
\frac{1}{p\left(  p+1\right)  } & \frac{1}{p\left(  p+1\right)  } & \frac
{1}{2}%
\end{array}
\right]
\]
with determinant%
\[
\det\mathbf{G}=\frac{\left(  p+p^{2}+4\right)  \left(  p+2\right)  ^{2}\left(
p-1\right)  ^{2}}{8p^{3}\left(  p+1\right)  ^{3}}.
\]
Since $p\geq3$, the determinant is positive and $\mathbf{G}$ is regular.

Let (\ref{Assumptab}) be satisfied. The functions $1_{\mathbf{z}}$ and
$q_{p-1,\mathbf{z}^{\prime}}$, $\mathbf{z}^{\prime}\in\mathcal{C}_{\mathbf{z}%
}$, are linearly independent and belong to $N_{p,\mathbf{z}}$. If the
dimension formula (\ref{dimform}) holds, they form a basis of $N_{p,\mathbf{z}%
}$.%
%TCIMACRO{\TeXButton{End Proof}{\endproof}}%
%BeginExpansion
\endproof
%EndExpansion

\begin{lemma}
Let (\ref{Assumptab}) be satisfied. Then%
\begin{equation}
\left(  \operatorname*{span}\left\{  1_{\mathbf{z}}\right\}
+\operatorname*{span}\left\{  q_{p-1,\mathbf{z}^{\prime}},\mathbf{z}^{\prime
}\in\mathcal{C}_{\mathbf{z}}\right\}  \right)  \cap N_{p,\mathbf{z}%
}^{\operatorname*{CR}}=\operatorname*{span}\left\{  1_{\mathbf{z}}\right\}  .
\label{claimfirstmain}%
\end{equation}
If the dimension formula (\ref{dimform}) holds, then $\dim N_{p,\mathbf{z}%
}^{\operatorname*{CR}}=1$.
\end{lemma}

%

%TCIMACRO{\TeXButton{Proof}{\proof}}%
%BeginExpansion
\proof
%EndExpansion
We distinguish between different cases.

\textbf{Case 1: }$\mathbf{z}$ is a critical point in $\mathcal{T}_{\mathbf{z}%
}$.

In this case we have the geometric situation as described in case 3b of the
proof of Lemma \ref{Lembasiscontker}. In particular there is no further
critical point in $\mathcal{T}_{\mathbf{z}}$ and the left-hand side in
(\ref{claimfirstmain}) equals $\operatorname*{span}\left\{  1_{\mathbf{z}%
},q_{p-1,\mathbf{z}}\right\}  $. We choose $\eta=\alpha1_{\mathbf{z}}+\beta
q_{p-1,\mathbf{z}}$. Our goal is to prove the implication:%
\begin{equation}
\eta\in N_{p,\mathbf{z}}^{\operatorname*{CR}}\implies\beta=0. \label{kernimpl}%
\end{equation}
Let $E\in\mathcal{E}_{\mathbf{z}}$ with $\mathbf{z}\in E$ and let
$K,K^{\prime}$ be the adjacent triangles in $\mathcal{T}_{\mathbf{z}}$; see
Fig. \ref{Figcrisscross}%
%TCIMACRO{\FRAME{ftbpFU}{2.316in}{2.1015in}{0pt}{\Qcb{Geometric situation,
%where $\QTR{bf}{z}$ is a critical vertex in the patch
%$\QTR{cal}{T}_{\QTR{bf}{z}}$.}}{\Qlb{Figcrisscross}}{crisscrosssing.eps}%
%{\special{ language "Scientific Word";  type "GRAPHIC";
%maintain-aspect-ratio TRUE;  display "USEDEF";  valid_file "F";
%width 2.316in;  height 2.1015in;  depth 0pt;  original-width 2.2762in;
%original-height 2.0626in;  cropleft "0";  croptop "1";  cropright "1";
%cropbottom "0";  filename 'crisscrosssing.EPS';file-properties "XNPEU";}} }%
%BeginExpansion
%EndExpansion
. Let $\mathbf{A}_{i}$ and $\mathbf{A}_{i}^{\prime}$, $1\leq i\leq3$, denote
their vertices with the convention that $\mathbf{A}_{1}=\mathbf{A}_{1}%
^{\prime}=\mathbf{z}$. The angle in $K$ at $\mathbf{A}_{i}$ is denoted by
$\alpha_{i}$ and $\alpha_{i}^{\prime}$ is the angle in $K^{\prime}$ at
$\mathbf{A}_{i}^{\prime}$. The edge $E_{i}$ is opposite to the vertex
$\mathbf{A}_{i}$ and $E_{i}^{\prime}$ is opposite to $\mathbf{A}_{i}^{\prime}%
$. Hence, $E=E_{2}=E_{3}^{\prime}$. We fix the vector $\mathbf{v}_{E}$ in
(\ref{basisvelodd}) as the unit vector which is orthogonal to the edge $E$ and
points into $K^{\prime}$. This vector will be denoted by $\mathbf{n}_{E}$. In
view of (\ref{defNCR}), implication (\ref{kernimpl}) follows if we prove
\begin{equation}
\left(  \eta,\operatorname*{div}\left(  B_{p,E}^{\operatorname*{CR}}%
\mathbf{n}_{E}\right)  \right)  _{L^{2}\left(  \omega_{\mathbf{z}}\right)
}=0\implies\beta=0. \label{kernimpl2}%
\end{equation}
By an integration by parts we get%
\[
\left(  \eta,\operatorname*{div}\left(  B_{p,E}^{\operatorname*{CR}}%
\mathbf{n}_{E}\right)  \right)  _{L^{2}\left(  \omega_{\mathbf{z}}\right)
}=\beta\left(  q_{p-1,\mathbf{z}},\operatorname*{div}\left(  B_{p,E}%
^{\operatorname*{CR}}\mathbf{n}_{E}\right)  \right)  _{L^{2}\left(  \omega
_{E}\right)  }%
\]
and (\ref{kernimpl2}) follows if we prove
\begin{equation}
\left(  q_{p-1,\mathbf{z}},\operatorname*{div}\left(  B_{p,E}%
^{\operatorname*{CR}}\mathbf{n}_{E}\right)  \right)  _{L^{2}\left(  \omega
_{E}\right)  }\neq0. \label{kernimpl3}%
\end{equation}
We get%
\begin{align*}
\left(  q_{p-1,\mathbf{z}},\operatorname*{div}\left(  B_{p,E}%
^{\operatorname*{CR}}\mathbf{n}_{E}\right)  \right)  _{L^{2}\left(  \omega
_{E}\right)  }  &  =\frac{\sigma_{K}}{\left\vert K\right\vert }\int_{K}%
P_{p-1}^{\left(  0,2\right)  }\left(  1-2\lambda_{K,1}\right)  \partial
_{\mathbf{n}_{E}}B_{p,E}^{\operatorname*{CR}}\\
&  +\frac{\sigma_{K^{\prime}}}{\left\vert K^{\prime}\right\vert }%
\int_{K^{\prime}}P_{p-1}^{\left(  0,2\right)  }\left(  1-2\lambda_{K^{\prime
},1}\right)  \partial_{\mathbf{n}_{E}}B_{p,E}^{\operatorname*{CR}}.
\end{align*}

For the first integral, we get explicitly%
\begin{equation}
\int_{K}P_{p-1}^{\left(  0,2\right)  }\left(  1-2\lambda_{K,1}\right)
\partial_{\mathbf{n}_{E}}B_{p,E}^{\operatorname*{CR}}=-2\partial
_{\mathbf{n}_{E}}\lambda_{K,2}\int_{K}P_{p-1}^{\left(  0,2\right)  }\left(
1-2\lambda_{K,1}\right)  \left(  P_{p}^{\left(  0,0\right)  }\right)
^{\prime}\left(  1-2\lambda_{K,2}\right)  . \label{Jacobimixed}%
\end{equation}
We employ the affine map $\chi:K\rightarrow K^{\prime}$ characterized by
$\chi\left(  \mathbf{A}_{1}\right)  =\mathbf{A}_{1}^{\prime}$, $\chi\left(
\mathbf{A}_{3}\right)  =\mathbf{A}_{2}^{\prime}$, $\chi\left(  \mathbf{A}%
_{2}\right)  =\mathbf{A}_{3}^{\prime}$. By taking into account $\sigma
_{K^{\prime}}=-\sigma_{K}$ we get%
\begin{align}
&  \left(  q_{p-1,\mathbf{z}},\operatorname*{div}\left(  B_{p,E}%
^{\operatorname*{CR}}\mathbf{n}_{E}\right)  \right)  _{L^{2}\left(  \omega
_{E}\right)  }=-2\frac{\sigma_{K}}{\left\vert K\right\vert }\left(
\partial_{\mathbf{n}_{E}}\lambda_{K,2}-\partial_{\mathbf{n}_{E}}%
\lambda_{K^{\prime},3}\right)  \iota_{K,p}\label{CRcritpointcent}\\
&  \text{with }\iota_{K,p}:=\int_{K}P_{p-1}^{\left(  0,2\right)  }\left(
1-2\lambda_{K,1}\right)  \left(  P_{p}^{\left(  0,0\right)  }\right)
^{\prime}\left(  1-2\lambda_{K,2}\right)  . \label{CRcritpointcent2}%
\end{align}
From (\ref{normcomp}) and elementary triangle trigonometry we conclude that%
\begin{align*}
&  \partial_{\mathbf{n}_{E}}\lambda_{K,2}-\partial_{\mathbf{n}_{E}}%
\lambda_{K^{\prime},3}=\partial_{\mathbf{n}_{2}}\lambda_{K,2}+\partial
_{\mathbf{n}_{2}}\lambda_{K^{\prime},3}=-\frac{\left\vert E\right\vert
}{2\left\vert K\right\vert }-\frac{\left\vert E\right\vert }{2\left\vert
K^{\prime}\right\vert }\\
&  \qquad=-\frac{1}{\sin\alpha_{1}\left\vert E_{3}\right\vert }-\frac{1}%
{\sin\alpha_{1}^{\prime}\left\vert E_{2}^{\prime}\right\vert }=-\frac
{\sin\alpha_{2}}{\left\vert E\right\vert \sin\alpha_{1}\sin\alpha_{3}}%
-\frac{\sin\alpha_{3}^{\prime}}{\left\vert E\right\vert \sin\alpha_{1}%
^{\prime}\sin\alpha_{2}^{\prime}}\\
&  \qquad=-\frac{1}{\left\vert E\right\vert }\left(  \frac{\sin\left(
\alpha_{1}+\alpha_{3}\right)  }{\sin\alpha_{1}\sin\alpha_{3}}+\frac
{\sin\left(  \alpha_{1}^{\prime}+\alpha_{2}^{\prime}\right)  }{\sin\alpha
_{1}^{\prime}\sin\alpha_{2}^{\prime}}\right)  =-\frac{1}{\left\vert
E\right\vert }\left(  \cot\alpha_{1}+\cot\alpha_{3}+\cot\alpha_{1}^{\prime
}+\cot\alpha_{2}^{\prime}\right)  .
\end{align*}
Recall that $\alpha_{1}+\alpha_{1}^{\prime}=\pi$ so that $\cot\alpha_{1}%
+\cot\alpha_{1}^{\prime}=0$ and%
\begin{equation}
\partial_{\mathbf{n}_{E}}\lambda_{K,2}-\partial_{\mathbf{n}_{E}}%
\lambda_{K^{\prime},3}=\partial_{\mathbf{n}_{2}}\lambda_{K,2}+\partial
_{\mathbf{n}_{2}}\lambda_{K^{\prime},3}=-\frac{1}{\left\vert E\right\vert
}\left(  \cot\alpha_{3}+\cot\alpha_{2}^{\prime}\right)  \neq0.
\label{normaljump1}%
\end{equation}

It remains to prove that the integral in (\ref{CRcritpointcent2}) is different
from zero. We obtain%
\begin{align*}
\iota_{K,p}  &  =2\left\vert K\right\vert \int_{0}^{1}\int_{0}^{1-x_{1}%
}P_{p-1}^{\left(  0,2\right)  }\left(  1-2x_{1}\right)  \left(  P_{p}^{\left(
0,0\right)  }\right)  ^{\prime}\left(  1-2x_{2}\right)  dx_{2}dx_{1}\\
&  =-\left\vert K\right\vert \int_{0}^{1}P_{p-1}^{\left(  0,2\right)  }\left(
1-2x_{1}\right)  \left(  \left.  P_{p}^{\left(  0,0\right)  }\left(
1-2x_{2}\right)  \right\vert _{0}^{1-x_{1}}\right)  dx_{1}\\
&  =-\frac{\left\vert K\right\vert }{2}\int_{-1}^{1}P_{p-1}^{\left(
0,2\right)  }\left(  t\right)  \left(  P_{p}^{\left(  0,0\right)  }\left(
-t\right)  -P_{p}^{\left(  0,0\right)  }\left(  1\right)  \right)  dt\\
&  =\frac{\left\vert K\right\vert }{2}\int_{-1}^{1}P_{p-1}^{\left(
0,2\right)  }\left(  t\right)  \left(  P_{p}^{\left(  0,0\right)  }\left(
t\right)  +1\right)  dt,
\end{align*}
where we used that $P_{p}^{\left(  0,0\right)  }$ is an odd function and
$P_{p}^{\left(  0,0\right)  }\left(  1\right)  =1$. The orthogonality of the
Legendre polynomial $P_{p}^{\left(  0,0\right)  }$ implies $\int_{-1}^{1}%
P_{p}^{\left(  0,0\right)  }\left(  t\right)  P_{p-1}^{\left(  0,2\right)
}\left(  t\right)  =0$ and from (\ref{02o}) we conclude that $\iota
_{K,p}=\left(  -1\right)  ^{p-1}\left\vert K\right\vert $. The combination
with (\ref{normaljump1}) leads to%
\begin{equation}
\left(  q_{p-1,\mathbf{z}^{\prime}},\operatorname*{div}\left(  B_{p,E}%
^{\operatorname*{CR}}\mathbf{n}_{E}\right)  \right)  _{L^{2}\left(
\omega_{\mathbf{z}}\right)  }=2\sigma_{K}\left(  -1\right)  ^{p-1}\left(
\frac{\cot\alpha_{3}+\cot\alpha_{2}^{\prime}}{\left\vert E\right\vert
}\right)  \neq0. \label{qkcommon}%
\end{equation}
In turn, the implication (\ref{kernimpl3}) follows and concludes the proof of
the assertion for the case that $\mathbf{z}$ is a critical point.

\textbf{Case 2. }$\mathbf{z}$ is not a critical point in $\mathcal{T}%
_{\mathbf{z}}$.

Let $\mathcal{C}_{\mathbf{z}}$ denote the set of critical points in
$\mathcal{T}_{\mathbf{z}}$. If $\mathcal{C}_{\mathbf{z}}$ is empty, we
conclude from (\ref{Lembasiscontker}), (\ref{inclN}), and (\ref{dimcond}) that
$\dim N_{p,\mathbf{z}}^{\operatorname*{CR}}=1$ and it remains to consider the
case $\mathcal{C}_{\mathbf{z}}\neq\emptyset$.

\textbf{Part a: Derivation of sufficient conditions.}

In view of (\ref{claimfirstmain}), we consider functions $\eta$ of the form:%
\[
\eta=\beta_{1}1_{\mathbf{z}}+\sum_{\mathbf{z}^{\prime}\in\mathcal{C}%
_{\mathbf{z}}}\beta_{\mathbf{z}^{\prime}}q_{p-1,\mathbf{z}^{\prime}}%
\quad\text{for real coefficients }\beta_{1}\text{, }\beta_{\mathbf{z}^{\prime
}}.
\]
Our goal is to prove the implication%
\[
\eta\in N_{p,\mathbf{z}}^{\operatorname*{CR}}\implies\forall\mathbf{z}%
^{\prime}\in\mathcal{C}_{\mathbf{z}}\quad\beta_{\mathbf{z}^{\prime}}=0.
\]
A sufficient condition is to prove the implication%
\[
\left(  \forall E\in\mathcal{E}_{\mathbf{z}}\text{ with }\mathbf{z}\in
E\mathbf{:\quad}\left(  \eta,\operatorname*{div}\left(  B_{p,E}%
^{\operatorname*{CR}}\mathbf{n}_{E}\right)  \right)  _{L^{2}\left(
\omega_{\mathbf{z}}\right)  }=0\right)  \implies\left(  \forall\mathbf{z}%
^{\prime}\in\mathcal{C}_{\mathbf{z}}:\quad\beta_{\mathbf{z}^{\prime}%
}=0\right)  .
\]
Since $\left(  1_{\mathbf{z}},\operatorname*{div}\left(  B_{p,E}%
^{\operatorname*{CR}}\mathbf{n}_{E}\right)  \right)  _{L^{2}\left(
\omega_{\mathbf{z}}\right)  }=0$ this is equivalent to showing%
\begin{align}
&  \left(  \forall E\in\mathcal{E}_{\mathbf{z}}\text{ with }\mathbf{z}\in
E\mathbf{:\quad}\sum_{\mathbf{z}^{\prime}\in\mathcal{C}_{\mathbf{z}}}%
\beta_{\mathbf{z}^{\prime}}\left(  q_{p-1,\mathbf{z}^{\prime}}%
,\operatorname*{div}\left(  B_{p,E}^{\operatorname*{CR}}\mathbf{n}_{E}\right)
\right)  _{L^{2}\left(  \omega_{\mathbf{z}}\right)  }=0\right)
\label{finalcondcrpo}\\
&  \implies\left(  \forall\mathbf{z}^{\prime}\in\mathcal{C}_{\mathbf{z}}%
:\quad\beta_{\mathbf{z}^{\prime}}=0\right)  .\nonumber
\end{align}

Next, we consider (\ref{finalcondcrpo}) for a reduced number of edges. Since
not all points in $\mathcal{V}_{\mathbf{z}}\backslash\left\{  \mathbf{z}%
\right\}  $ can be critical points, we can group the points in $\mathcal{C}%
_{\mathbf{z}}$ as follows (see Fig. \ref{Fig_nod_patch}).%
%TCIMACRO{\FRAME{ftbpFU}{3.0355in}{1.855in}{0pt}{\Qcb{Nodal patch, illustrating
%singular vertices at $\partial\omega_{\QTR{bf}{z}}$. They are partitioned into
%edge-connected parts $\QTR{cal}{C}_{\QTR{bf}{z},\ell}$, $1\leq\ell\leq
%c_{\QTR{bf}{z}}$. The extremal points $\QTR{bf}{z}_{\ell,0}$ and
%$\QTR{bf}{z}_{\ell,n_{\ell}+1}$ are no singular points. The edge connecting
%$\QTR{bf}{z}$ with $\QTR{bf}{z}_{\ell,j}$ is denoted by $E_{\ell,j}$ and the
%connecting edge for the singular points $\QTR{bf}{z}_{\ell,j-1}$,
%$\QTR{bf}{z}_{\ell,j}$ is denoted by $E_{\ell,j}^{0}$. The normal vector
%$\QTR{bf}{n}_{\ell,j}$ at $E_{\ell,j}$ is counterclockwise oriented. Condition
%(\ref{finalcondcrpo}) is considered for all $\QTR{cal}{C}_{\QTR{bf}{z},\ell}$
%separately and the index $\ell$ is skipped in the computation.}}%
%{\Qlb{Fig_nod_patch}}{faecher_29032021.eps}%
%{\special{ language "Scientific Word";  type "GRAPHIC";
%maintain-aspect-ratio TRUE;  display "USEDEF";  valid_file "F";
%width 3.0355in;  height 1.855in;  depth 0pt;  original-width 2.0531in;
%original-height 1.2427in;  cropleft "0";  croptop "1";  cropright "1";
%cropbottom "0";  filename 'faecher_29032021.eps';file-properties "XNPEU";}} }%
%BeginExpansion
%EndExpansion
We say two critical points $\mathbf{z}^{\prime},\mathbf{z}^{\prime\prime}%
\in\mathcal{C}_{\mathbf{z}}$ are \textit{edge-connected} if there is an edge
$E\in\mathcal{E}_{\mathbf{z}}$ having $\mathbf{z}^{\prime},\mathbf{z}%
^{\prime\prime}$ as endpoints. We decompose $\mathcal{C}_{\mathbf{z}}$ into a
minimal number of subsets $\mathcal{C}_{\mathbf{z},\ell}$, $1\leq\ell\leq
c_{\mathbf{z}}$, such that all points in $\mathcal{C}_{\mathbf{z},\ell}$ are
edge-connected in that the points in $\mathcal{C}_{\mathbf{z},\ell}$ can be
numbered counterclockwise $\mathbf{z}_{\ell,1},\mathbf{z}_{\ell,2}%
,\ldots,\mathbf{z}_{\ell,n_{\ell}}$ and satisfy: $\mathbf{z}_{\ell,j-1}$ and
$\mathbf{z}_{\ell,j}$ are edge-connected for $j=2,3,\ldots n_{\ell}$. The
construction implies that $\mathbf{z}_{\ell,1}$ is edge-connected with some
$\mathbf{z}_{\ell,0}\in\mathcal{V}_{\mathbf{z}}\backslash\left\{
\mathbf{z}\right\}  $, which is not a critical point, and $\mathbf{z}%
_{\ell,n_{\ell}}$ is edge-connected to some point $\mathbf{z}_{\ell,n_{\ell
}+1}\in\mathcal{V}_{\mathbf{z}}\backslash\left\{  \mathbf{z}\right\}  $, which
is also no critical point. In (\ref{finalcondcrpo}), we choose the edges
$E=E_{\ell,j}$, which is edge-connected to $\mathbf{z}_{\ell,j}$ with
$\mathbf{z}$, $1\leq j\leq n_{\ell}$. The edge connecting $\mathbf{z}%
_{\ell,j-1}$ with $\mathbf{z}_{\ell,j}$ is denoted by $E_{\ell,j}^{0}$ and the
triangle with vertices $\mathbf{z}$, $\mathbf{z}_{\ell,j-1}$, $\mathbf{z}%
_{\ell,j}$ by $K_{\ell,j}$. Let $\mathbf{n}_{\ell,j}$ denote the unit vector
which is orthogonal to $E_{\ell,j}$ and points into $K_{\ell,j+1}$. The angle
in $K_{\ell,j}$ at the vertex $\mathbf{z}$ is locally denoted by $\alpha
_{\ell,1}$ and at the vertices $\mathbf{z}_{\ell,j-1}$, $\mathbf{z}_{\ell,j}$
by $\alpha_{\ell,2}$, $\alpha_{\ell,3}$.

We write short $B_{\ell,j}^{\operatorname*{CR}}$ for the test function
$B_{p,E}^{\operatorname*{CR}}$ with $E=E_{\ell,j}$. For fixed $\ell$, we test
with the functions $B_{\ell,j}^{\operatorname*{CR}}\mathbf{n}_{\ell,j}$,
$1\leq j\leq n_{\ell}$, and obtain the conditions%
\begin{equation}
\sum_{j\mathbf{=}1}^{n_{\ell}}\beta_{\ell,j}\left(  q_{p-1,\mathbf{z}_{\ell
,j}},\operatorname*{div}\left(  B_{\ell,i}^{\operatorname*{CR}}\mathbf{n}%
_{\ell,j}\right)  \right)  _{L^{2}\left(  \omega_{\mathbf{z}}\right)  }%
=0\quad\forall1\leq i\leq n_{\ell}. \label{betacond}%
\end{equation}
Hence, (\ref{finalcondcrpo}) is proved if we show the sufficient condition%
\begin{equation}
\forall1\leq\ell\leq c_{\ell}\quad\text{(\ref{betacond}) has the unique
solution }\beta_{\ell,i}=0\quad\forall1\leq i\leq n_{\ell}. \label{impl4}%
\end{equation}
Since these conditions are decoupled with respect $\ell$ we verify this for
each $\ell$ separately and drop the index $\ell$. We define the matrix
$\mathbf{M}_{n}=\left(  m_{i,j}\right)  _{\substack{1\leq i\leq n\\1\leq j\leq
n}}$ by
\[
m_{i,j}:=\left(  q_{p-1,\mathbf{z}_{j}},\operatorname*{div}\left(
B_{i}^{\operatorname*{CR}}\mathbf{n}_{i}\right)  \right)  _{L^{2}\left(
\omega_{\mathbf{z}}\right)  }.
\]
By using a) the support properties of the critical functions
$q_{p-1,\mathbf{z}_{j}}$ and the Crouzeix-Raviart functions $B_{\ell
,i}^{\operatorname*{CR}}$, b) the fixed counterclockwise orientation of the
normal vector $\mathbf{n}_{i}$, c) the fact that $B_{i}^{\operatorname*{CR}}$
and $q_{p-1,\mathbf{z}_{j}}$ are defined by barycentric coordinates, d) the
function $\sigma_{K}$ has alternating sign, we obtain that $\mathbf{M}_{n}$
has the following tridiagonal structure%
\begin{equation}
\mathbf{M}_{n}=\left[
\begin{array}
[c]{ccccc}%
d_{1} & d_{12} & 0 & \ldots & 0\\
d_{21} & d_{2} & d_{23} & \ddots & \vdots\\
0 & d_{32} & d_{3} & \ddots & 0\\
\vdots & \ddots & \ddots & \ddots & d_{n-1,n}\\
0 & \ldots & 0 & d_{n,n-1} & d_{n}%
\end{array}
\right]  \label{defMn}%
\end{equation}
and condition (\ref{impl4}) is equivalent to showing that $\mathbf{M}_{n}$ is
regular. In the following we compute the matrix entries of $\mathbf{M}_{n}%
$.\medskip

\textbf{Part b. Computing the diagonal entries }$d_{j}$.

Now, let $\mathbf{z}^{\prime}\in\mathcal{C}_{\mathbf{z}}$ and denote by $E$
the edge connecting $\mathbf{z}$ with $\mathbf{z}^{\prime}$. The adjacent
triangles are denoted by $K,K^{\prime}$ with vertices $\mathbf{A}_{i}$,
$\mathbf{A}_{i}^{\prime}$, $i=1,2,3$, $\mathbf{A}_{1}=\mathbf{A}_{1}^{\prime
}=\mathbf{z}$, and $\mathbf{A}_{3}=\mathbf{A}_{2}^{\prime}=\mathbf{z}^{\prime
}$. We first evaluate%
\[
\left(  q_{p-1,\mathbf{z}^{\prime}},\operatorname*{div}\left(  B_{p,E}%
^{\operatorname*{CR}}\mathbf{n}_{E}\right)  \right)  _{L^{2}\left(
\omega_{\mathbf{z}}\right)  }.
\]
By repeating the computation as in case 1 we conclude as for (\ref{qkcommon})
that this quantity satisfies%
\[
\left(  q_{p-1,\mathbf{z}^{\prime}},\operatorname*{div}\left(  B_{p,E}%
^{\operatorname*{CR}}\mathbf{n}_{E}\right)  \right)  _{L^{2}\left(
\omega_{\mathbf{z}}\right)  }=2\sigma_{K}\left(  -1\right)  ^{p-1}\frac
{\cot\alpha_{1}+\cot\alpha_{1}^{\prime}}{\left\vert E\right\vert }.
\]

Let $K_{1}\in\mathcal{T}_{\mathbf{z}}$ denote the triangle with vertices
$\mathbf{z}$, $\mathbf{z}_{0}$, $\mathbf{z}_{1}$, i.e., $\mathbf{z}_{0}$ is
not a critical point. Without loss of generality we assume that the signature
of $K_{1}$ satisfies $\sigma_{K_{1}}=-1$ while the proof for $\sigma_{K_{1}%
}=1$ is verbatim. By this signature convention and since $p$ is odd, it
follows that%
\[
d_{j}=2\left(  -1\right)  ^{j}\frac{\cot\alpha_{j,1}+\cot\alpha_{j+1,1}%
}{\left\vert E_{j}\right\vert }.
\]

\textbf{Part c. Computing the off-diagonal entries}.

Let $\mathbf{z}^{\prime}\in\mathcal{E}_{\mathbf{z}}$ be a critical point as
before and recall the notation of $E$, $\mathbf{A}_{i}$, $\mathbf{A}%
_{i}^{\prime}$ as in part b, in particular, $\mathbf{A}_{1}=\mathbf{A}%
_{1}^{\prime}=\mathbf{z}$, $\mathbf{A}_{3}=\mathbf{A}_{2}^{\prime}%
=\mathbf{z}^{\prime}$, $K$ is the triangle with vertices $\mathbf{A}_{i}$,
$1\leq i\leq3$, and $E$ has endpoints $\mathbf{z}$ and $\mathbf{z}^{\prime}$.
We assume that $\mathbf{A}_{2}\in\mathcal{C}_{\mathbf{z}}$ so that $\left(
q_{p-1,\mathbf{A}_{2}},\operatorname*{div}\left(  B_{p,E}^{\operatorname*{CR}%
}\mathbf{n}_{E}\right)  \right)  _{L^{2}\left(  \omega_{\mathbf{z}}\right)  }$
is a matrix entry of the form $d_{j,j-1}$. Then, we evaluate%
\begin{align*}
&  \left(  q_{p-1,\mathbf{A}_{2}},\operatorname*{div}\left(  B_{p,E}%
^{\operatorname*{CR}}\mathbf{n}_{E}\right)  \right)  _{L^{2}\left(
\omega_{\mathbf{z}}\right)  }=\left(  q_{p-1,\mathbf{A}_{2}}%
,\operatorname*{div}\left(  B_{p,E}^{\operatorname*{CR}}\mathbf{n}_{E}\right)
\right)  _{L^{2}\left(  K\right)  }\\
&  \qquad=-2\sigma_{K}\left(  \partial_{\mathbf{n}_{E}}\lambda_{K,\mathbf{A}%
_{2}}\right)  \theta_{p}.
\end{align*}
with%
\[
\theta_{p}:=\frac{1}{\left\vert K\right\vert }\int_{K}P_{p-1}^{\left(
0,2\right)  }\left(  1-2\lambda_{K,\mathbf{A}_{2}}\right)  \left(
P_{p}^{\left(  0,0\right)  }\right)  ^{\prime}\left(  1-2\lambda
_{K,\mathbf{A}_{2}}\right)  .
\]
The integral in the definition of $\theta_{p}$ is transformed to the reference
element by the affine transform $\chi$ characterized by%
\begin{equation}
\chi\left(  \mathbf{\hat{A}}_{1}\right)  =\mathbf{A}_{1}=\mathbf{z},\quad
\chi\left(  \mathbf{\hat{A}}_{2}\right)  =\mathbf{A}_{2},\quad\chi\left(
\mathbf{\hat{A}}_{3}\right)  =\mathbf{A}_{3}. \label{defchiloc1}%
\end{equation}
This implies in particular that $\theta_{p}$ is independent of $K$. Thus,%
\[
\left(  q_{p-1,\mathbf{A}_{2}},\operatorname*{div}\left(  B_{p,E}%
^{\operatorname*{CR}}\mathbf{n}_{E}\right)  \right)  _{L^{2}\left(
\omega_{\mathbf{z}}\right)  }=-2\sigma_{K}\theta_{p}\partial_{\mathbf{n}_{E}%
}\lambda_{K,\mathbf{A}_{2}}.
\]
From Appendix \ref{DerBary} (see (\ref{normcomp})) it follows that%
\[
\partial_{\mathbf{n}_{E}}\lambda_{K,\mathbf{A}_{2}}=-\frac{\left\vert
E\right\vert }{2\left\vert K\right\vert }%
\]
so that%
\[
\left(  q_{p-1,\mathbf{A}_{2}},\operatorname*{div}\left(  B_{p,E}%
^{\operatorname*{CR}}\mathbf{n}_{E}\right)  \right)  _{L^{2}\left(
\omega_{\mathbf{z}}\right)  }=\frac{\sigma_{K}\left\vert E\right\vert
}{\left\vert K\right\vert }\theta_{p}.
\]
By using the affine transform $\chi$ as in (\ref{defchiloc1}) we are led to%
\begin{align*}
\theta_{p}  &  =2\int_{0}^{1}\int_{0}^{1-x_{1}}P_{p-1}^{\left(  0,2\right)
}\left(  1-2x_{1}\right)  \left(  P_{p}^{\left(  0,0\right)  }\right)
^{\prime}\left(  1-2x_{1}\right)  dx_{2}dx_{1}\\
&  =2\int_{0}^{1}\left(  1-x_{1}\right)  P_{p-1}^{\left(  0,2\right)  }\left(
1-2x_{1}\right)  \left(  P_{p}^{\left(  0,0\right)  }\right)  ^{\prime}\left(
1-2x_{1}\right)  dx_{1}\\
&  =\frac{1}{2}\int_{-1}^{1}\left(  t+1\right)  P_{p-1}^{\left(  0,2\right)
}\left(  t\right)  \left(  P_{p}^{\left(  0,0\right)  }\right)  ^{\prime
}\left(  t\right)  dt\\
&  \overset{\text{\cite[Tab. 18.6.1]{NIST:DLMF}}}{=}\frac{1}{2}\int_{-1}%
^{1}\left(  1-t\right)  P_{p-1}^{\left(  2,0\right)  }\left(  t\right)
\left(  P_{p}^{\left(  0,0\right)  }\right)  ^{\prime}\left(  t\right)  dt\\
&  \overset{\text{\cite[18.9.15]{NIST:DLMF}}}{=}\dfrac{p+1}{4}\int_{-1}%
^{1}\left(  1-t\right)  P_{p-1}^{\left(  2,0\right)  }\left(  t\right)
P_{p-1}^{(1,1)}\left(  t\right)  dt\\
&  \overset{\text{(\ref{211b})}}{=}1.
\end{align*}
Hence,%
\[
\left(  q_{p-1,\mathbf{A}_{2}},\operatorname*{div}\left(  B_{p,E}%
^{\operatorname*{CR}}\mathbf{n}_{E}\right)  \right)  _{L^{2}\left(
\omega_{\mathbf{z}}\right)  }=\frac{\sigma_{K}\left\vert E\right\vert
}{\left\vert K\right\vert }%
\]
and we conclude that%
\[
d_{i,i-1}=\left(  -1\right)  ^{i}\frac{\left\vert E_{i}\right\vert
}{\left\vert K_{i}\right\vert }%
\]
holds. By employing an affine transform $\chi$, which maps $K_{i+1}$ onto
$K_{i}$ and is the identity on $E_{i}$, we can reuse the computations for
$d_{i,i-1}$ to get%
\[
d_{i,i+1}=\left(  -1\right)  ^{i}\frac{\left\vert E_{i}\right\vert
}{\left\vert K_{i+1}\right\vert }.
\]
In summary,%
\begin{align*}
&  \mathbf{M}_{n}=\\
&  =2\left[
\begin{array}
[c]{ccccc}%
-\frac{\cot\alpha_{1,1}+\cot\alpha_{2,1}}{\left\vert E_{1}\right\vert } &
-\frac{\left\vert E_{1}\right\vert }{2\left\vert K_{2}\right\vert } & 0 &
\ldots & 0\\
\frac{\left\vert E_{2}\right\vert }{2\left\vert K_{2}\right\vert } &
\frac{\cot\alpha_{2,1}+\cot\alpha_{3,1}}{\left\vert E_{2}\right\vert } &
\frac{\left\vert E_{2}\right\vert }{2\left\vert K_{3}\right\vert } & \ddots &
\vdots\\
0 & -\frac{\left\vert E_{3}\right\vert }{2\left\vert K_{3}\right\vert } &
-\frac{\cot\alpha_{3,1}+\cot\alpha_{4,1}}{\left\vert E_{3}\right\vert } &
\ddots & 0\\
\vdots & \ddots & \ddots & \ddots & \frac{\left(  -1\right)  ^{n-1}\left\vert
E_{n-1}\right\vert }{2\left\vert K_{n}\right\vert }\\
0 & \ldots & 0 & \frac{\left(  -1\right)  ^{n}\left\vert E_{n}\right\vert
}{2\left\vert K_{n}\right\vert } & \left(  -1\right)  ^{n}\frac{\cot
\alpha_{n,1}+\cot\alpha_{n+1,1}}{\left\vert E_{n}\right\vert }%
\end{array}
\right] \\
&  =2\mathbf{D}_{\sigma}\mathbf{T}_{n,%
%TCIMACRO{\TeXButton{boldalpha}{\mbox{\boldmath$ \alpha$}}}%
%BeginExpansion
\mbox{\boldmath$ \alpha$}%
%EndExpansion
}\mathbf{D}_{\mathcal{E}}.
\end{align*}
with the diagonal matrices%
\[
\mathbf{D}_{\sigma}:=\operatorname*{diag}\left[  \left(  -1\right)  ^{i}:1\leq
i\leq n\right]  \quad\text{and\quad}\mathbf{D}_{\mathcal{E}}%
:=\operatorname*{diag}\left[  \left\vert E_{i}\right\vert ^{-1}:1\leq i\leq
n\right]
\]
and the symmetric tridiagonal matrix%
\begin{equation}
\mathbf{T}_{n,%
%TCIMACRO{\TeXButton{boldalpha}{\mbox{\boldmath$ \alpha$}}}%
%BeginExpansion
\mbox{\boldmath$ \alpha$}%
%EndExpansion
}:=\left[
\begin{array}
[c]{ccccc}%
\frac{\sin\left(  \alpha_{1,1}+\alpha_{2,1}\right)  }{\sin\alpha_{1,1}%
\sin\alpha_{2,1}} & \frac{1}{\sin\alpha_{2,1}} & 0 & \ldots & 0\\
\frac{1}{\sin\alpha_{2,1}} & \frac{\sin\left(  \alpha_{2,1}+\alpha
_{3,1}\right)  }{\sin\alpha_{2,1}\sin\alpha_{3,1}} & \frac{1}{\sin\alpha
_{3,1}} & \ddots & \vdots\\
0 & \frac{1}{\sin\alpha_{3,1}} & \frac{\sin\left(  \alpha_{3,1}+\alpha
_{4,1}\right)  }{\sin\alpha_{3,1}\sin\alpha_{4,1}} & \ddots & 0\\
\vdots & \ddots & \ddots & \ddots & \frac{1}{\sin\alpha_{n,1}}\\
0 & \ldots & 0 & \frac{1}{\sin\alpha_{n,1}} & \frac{\sin\left(  \alpha
_{n,1}+\alpha_{n+1,1}\right)  }{\sin\alpha_{n,1}\sin\alpha_{n+1,1}}%
\end{array}
\right]  . \label{defTnalpha}%
\end{equation}
By induction (see Appendix \ref{SecMatDet}) one proves%
\begin{equation}
\det\mathbf{T}_{n,%
%TCIMACRO{\TeXButton{boldalpha}{\mbox{\boldmath$ \alpha$}}}%
%BeginExpansion
\mbox{\boldmath$ \alpha$}%
%EndExpansion
}=\frac{\sin\left(  \sum_{i=1}^{n+1}\alpha_{i,1}\right)  }{\Pi_{i=1}^{n+1}%
\sin\alpha_{i,1}}. \label{detclaim}%
\end{equation}
We know that $0<\sum_{i=1}^{n+1}\alpha_{i,1}<\pi$ so that $\det\mathbf{T}_{n,%
%TCIMACRO{\TeXButton{boldalpha}{\mbox{\boldmath$ \alpha$}}}%
%BeginExpansion
\mbox{\boldmath$ \alpha$}%
%EndExpansion
}\mathbf{\neq0}$ and, in turn, $\det\mathbf{M}_{n}\neq0$.

Since $\mathbf{M}_{n}$ is the system matrix in the linear equations
(\ref{impl4}) for a fixed but arbitrary $\ell$ in (\ref{impl4}), we have
proved the implication (\ref{impl4}) and, in turn, the following statement: If
the dimension formula holds, we conclude from Lemma \ref{Lembasiscontker} that%
\[
N_{p,\mathbf{z}}=\operatorname*{span}\left\{  1_{\mathbf{z}}\right\}
+\operatorname*{span}\left\{  q_{p-1,\mathbf{z}^{\prime}}:\mathbf{z}^{\prime
}\in\mathcal{C}_{\mathbf{z}}\right\}  .
\]
Since $N_{p,\mathbf{z}}^{\operatorname*{CR}}\subset N_{p,\mathbf{z}}$ (cf.
(\ref{inclN})) we have in this case%
\[
N_{p,\mathbf{z}}^{\operatorname*{CR}}\subset\left(  \operatorname*{span}%
\left\{  1_{\mathbf{z}}\right\}  +\operatorname*{span}\left\{
q_{p-1,\mathbf{z}^{\prime}}:\mathbf{z}^{\prime}\in\mathcal{C}_{\mathbf{z}%
}\right\}  \right)  \cap N_{p,\mathbf{z}}^{\operatorname*{CR}}%
=\operatorname*{span}\left\{  1_{\mathbf{z}}\right\}
\]
and, in turn, $N_{p,\mathbf{z}}^{\operatorname*{CR}}=\operatorname*{span}%
\left\{  1_{\mathbf{z}}\right\}  $ (cf. Rem. \ref{Reminclusions}).%
%TCIMACRO{\TeXButton{End Proof}{\endproof}}%
%BeginExpansion
\endproof
%EndExpansion

\appendix

\section{Derivatives of Barycentric Coordinates\label{DerBary}}

As before we denote the vertices of a triangle $K$ by $\mathbf{A}_{i}$, $1\leq
i\leq3$, and set for ease of notation $\mathbf{A}_{0}:=\mathbf{A}_{3}$ and
$\mathbf{A}_{4}:=\mathbf{A}_{1}$. Then, for the barycentric coordinate
$\lambda_{K,\mathbf{A}_{i}}$ the following relations hold%

\begin{align}
\lambda_{K,\mathbf{A}_{i}}\left(  \mathbf{x}\right)   &  =\frac{\det\left[
\mathbf{x}-\mathbf{A}_{i-1}\mid\mathbf{A}_{i+1}-\mathbf{A}_{i-1}\right]
}{2\left\vert K\right\vert }\nonumber\\
\nabla\lambda_{K,\mathbf{A}_{i}}  &  =\nabla\frac{\left\langle \mathbf{x}%
-\mathbf{A}_{i-1},\left(  \mathbf{A}_{i+1}-\mathbf{A}_{i-1}\right)  ^{\perp
}\right\rangle }{2\left\vert K\right\vert }=\frac{\left(  \mathbf{A}%
_{i+1}-\mathbf{A}_{i-1}\right)  ^{\perp}}{2\left\vert K\right\vert }
\label{formulagrad}%
\end{align}
with $\mathbf{v}^{\perp}=\left(  v_{2},-v_{1}\right)  ^{T}$. The outer normal
vector $\mathbf{n}_{i}$ for the edge $E_{i}$ (opposite to $\mathbf{A}_{i}$) is
given by $\frac{\left(  \mathbf{A}_{i-1}-\mathbf{A}_{i+1}\right)  ^{\perp}%
}{\left\vert E_{i}\right\vert }.$ Hence,%
\begin{align}
\partial_{\mathbf{n}_{k}}\lambda_{K,\mathbf{A}_{i}}  &  =\frac{\left(
\mathbf{A}_{i+1}-\mathbf{A}_{i-1}\right)  ^{\perp}}{2\left\vert K\right\vert
}\frac{\left(  \mathbf{A}_{k-1}-\mathbf{A}_{k+1}\right)  ^{\perp}}{\left\vert
E_{k}\right\vert }=\frac{\left\langle \mathbf{A}_{i+1}-\mathbf{A}%
_{i-1},\mathbf{A}_{k-1}-\mathbf{A}_{k+1}\right\rangle }{2\left\vert
K\right\vert \left\vert E_{k}\right\vert }\label{normcomp}\\
&  =\frac{\cos\beta_{i,k}}{2\left\vert K\right\vert }\left\vert E_{i}%
\right\vert ,\nonumber
\end{align}
where $\beta_{i,k}\in\left]  0,\pi\right]  $ is the angle at that vertex in
$K$, where $E_{i}$ and $E_{k}$ meet with the convention $\beta_{i,i}=\pi$ for
all $i=1,2,3$. The tangential vector (counterclockwise orientation) for the
edge $i$ is given by $\mathbf{t}_{i}:=\frac{\mathbf{A}_{i-1}-\mathbf{A}_{i+1}%
}{\left\vert E_{i}\right\vert }$ so that%
\begin{align}
\partial_{\mathbf{t}_{k}}\lambda_{K,\mathbf{A}_{i}}  &  =\frac{\left(
\mathbf{A}_{i+1}-\mathbf{A}_{i-1}\right)  ^{\perp}}{2\left\vert K\right\vert
}\frac{\left(  \mathbf{A}_{k-1}-\mathbf{A}_{k+1}\right)  }{\left\vert
E_{k}\right\vert }\label{tangcomp}\\
&  =\frac{\det\left[  \mathbf{A}_{k-1}-\mathbf{A}_{k+1}\mid\mathbf{A}%
_{i+1}-\mathbf{A}_{i-1}\right]  }{2\left\vert K\right\vert \left\vert
E_{k}\right\vert }=\frac{\varepsilon_{i,k}}{\left\vert E_{k}\right\vert
}\nonumber
\end{align}
with%
\begin{equation}
\varepsilon_{i,k}:=\left(  k-i+1\right)  \operatorname{mod}3-1=\left\{
\begin{array}
[c]{ll}%
-1 & k=i-1,\\
0 & k=i,\\
1 & k=i+1.
\end{array}
\right.  \label{defepsik}%
\end{equation}

\section{The Determinant of $\mathbf{T}_{n,%
%TCIMACRO{\TeXButton{boldalpha}{\mbox{\boldmath$ \alpha$}}}%
%BeginExpansion
\mbox{\boldmath$ \alpha$}%
%EndExpansion
}$ (cf. (\ref{defTnalpha}))\label{SecMatDet}}

We write short $\mathbf{T}_{n}$ for $\mathbf{T}_{n,%
%TCIMACRO{\TeXButton{boldalpha}{\mbox{\boldmath$ \alpha$}}}%
%BeginExpansion
\mbox{\boldmath$ \alpha$}%
%EndExpansion
}$ and denote the matrix entries in $\mathbf{T}_{n}$ by $\left(
t_{n,i,j}\right)  _{i,j=1}^{n}$. We apply the well known three-term recursion
formula for the determinant of a symmetric tridiagonal matrix (see, e.g.,
\cite[(5.5.3)]{Stoer_Bulirsch_engl}) and obtain%
\[
\det\mathbf{T}_{n}=t_{n,n,n}\det\left(  \mathbf{T}_{n-1}\right)
-t_{n,n,n-1}^{2}\det\mathbf{T}_{n-2}.
\]
For $n=1,2$, it is straightforward to compute%
\[
\det\mathbf{T}_{n}=\frac{\sin\left(  \sum_{i=1}^{n+1}\alpha_{i,1}\right)
}{\Pi_{i=1}^{n+1}\sin\alpha_{i,1}}.
\]
We assume inductively (\ref{detclaim}) holds true for $\tilde{n}%
=1,2,\ldots,n-1$ and prove the assertion for $n$. We get%
\begin{align}
&  \frac{\sin\left(  \alpha_{n,1}+\alpha_{n+1,1}\right)  }{\sin\alpha
_{n,1}\sin\alpha_{n+1,1}}\det\left(  \mathbf{T}_{n-1}\right)  -\frac{1}%
{\sin^{2}\alpha_{n,1}}\det\mathbf{T}_{n-2}\nonumber\\
&  =\frac{\sin\left(  \alpha_{n,1}+\alpha_{n+1,1}\right)  }{\sin\alpha
_{n,1}\sin\alpha_{n+1,1}}\frac{\sin\left(  \sum_{i=1}^{n}\alpha_{i,1}\right)
}{\Pi_{i=1}^{n}\sin\alpha_{i,1}}-\frac{1}{\sin^{2}\alpha_{n,1}}\frac
{\sin\left(  \sum_{i=1}^{n-1}\alpha_{i,1}\right)  }{\Pi_{i=1}^{n-1}\sin
\alpha_{i,1}}\nonumber\\
&  =\frac{1}{\sin\alpha_{n,1}\Pi_{i=1}^{n+1}\sin\alpha_{i,1}}\left(
\sin\left(  \alpha_{n,1}+\alpha_{n+1,1}\right)  \sin\left(  \sum_{i=1}%
^{n}\alpha_{i,1}\right)  -\sin\alpha_{n+1,1}\sin\left(  \sum_{i=1}^{n-1}%
\alpha_{i,1}\right)  \right) \nonumber\\
&  =\frac{1}{\Pi_{i=1}^{n+1}\sin\alpha_{i,1}}\left(  \cos\alpha_{n+1,1}%
\sin\left(  \sum_{i=1}^{n}\alpha_{i,1}\right)  +\right. \nonumber\\
&  \left.  +\sin\alpha_{n+1,1}\frac{\cos\alpha_{n,1}\sin\left(  \sum_{i=1}%
^{n}\alpha_{i,1}\right)  -\sin\left(  \sum_{i=1}^{n}\alpha_{i,1}-\alpha
_{n,1}\right)  }{\sin\alpha_{n,1}}\right)  . \label{detformulaproof}%
\end{align}
The term in the numerator of the last term in the right-hand side in
(\ref{detformulaproof}) equals%
\[
\cos\alpha_{n,1}\sin\left(  \sum_{i=1}^{n}\alpha_{i,1}\right)  -\sin\left(
\sum_{i=1}^{n}\alpha_{i,1}-\alpha_{n,1}\right)  =\cos\left(  \sum_{i=1}%
^{n}\alpha_{i,1}\right)  \sin\left(  \alpha_{n,1}\right)  .
\]
Hence, it follows%
\begin{align*}
&  \frac{\sin\left(  \alpha_{n,1}+\alpha_{n+1,1}\right)  }{\sin\alpha
_{n,1}\sin\alpha_{n+1,1}}\det\left(  \mathbf{T}_{n-1}\right)  -\frac{1}%
{\sin^{2}\alpha_{n,1}}\det\mathbf{T}_{n-2}\\
&  =\frac{1}{\Pi_{i=1}^{n+1}\sin\alpha_{i,1}}\left(  \cos\alpha_{n+1,1}%
\sin\left(  \sum_{i=1}^{n}\alpha_{i,1}\right)  +\sin\alpha_{n+1,1}\cos\left(
\sum_{i=1}^{n}\alpha_{i,1}\right)  \right) \\
&  =\frac{\sin\left(  \sum_{i=1}^{n+1}\alpha_{i,1}\right)  }{\Pi_{i=1}%
^{n+1}\sin\alpha_{i,1}}%
\end{align*}
and, in turn, formula (\ref{detclaim}) for the determinant of $\mathbf{T}_{n}$.

\section{Closed Form of some Integrals involving Jacobi
Polynomials\label{SecIntJP}}

In this section, we derive some explicit formulae for integrals involving
Jacobi polynomials.

\begin{lemma}
\label{LemCHJP}It holds for all $n\in\mathbb{N}_{0}$%
\begin{align}
\int_{-1}^{1}\left(  1+t\right)  P_{n}^{\left(  0,2\right)  }\left(  t\right)
dt  &  =4\frac{\left(  -1\right)  ^{n}}{\left(  n+1\right)  \left(
n+2\right)  },\label{formel1IntP1}\\
\int_{-1}^{1}\left(  1+t\right)  \left(  P_{n}^{\left(  0,2\right)  }\left(
t\right)  \right)  ^{2}dt  &  =2.\nonumber
\end{align}

\end{lemma}

%

%TCIMACRO{\TeXButton{Proof}{\proof}}%
%BeginExpansion
\proof
%EndExpansion
\textbf{Part a.}

We set $I_{n}:=\int_{-1}^{1}\left(  1+t\right)  P_{n}^{\left(  0,2\right)
}\left(  t\right)  dt$. For $n=0$, we have $P_{0}^{\left(  0,2\right)
}\left(  t\right)  =1$ and $I_{0}=2$. For $n\geq1$, we employ \cite[18.9.5
with $\alpha=0$ and $\beta=1$ therein]{NIST:DLMF} in the form%
\begin{equation}
2\left(  n+1\right)  P_{n}^{(0,1)}=(n+2)P_{n}^{(0,2)}+nP_{n-1}^{(0,2)}.
\label{recJacobi}%
\end{equation}
The orthogonality relations of $P_{n}^{\left(  0,1\right)  }$ imply that
$\int_{-1}^{1}\left(  1+t\right)  P_{n}^{\left(  0,1\right)  }\left(
t\right)  dt=0$ so that the recurrence%

\[
0=\left(  n+2\right)  I_{n}+nI_{n-1}%
\]
follows. This can be resolved and we obtain%
\[
I_{n}=-\frac{n}{n+2}I_{n-1}=\left(  -1\right)  ^{n}2\frac{n!}{\left(
n+2\right)  !}I_{0}=4\frac{\left(  -1\right)  ^{n}}{\left(  n+1\right)
\left(  n+2\right)  }.
\]

\textbf{Part b}. We set $J_{n}:=\int_{-1}^{1}\left(  P_{n}^{\left(
0,2\right)  }\left(  t\right)  \right)  ^{2}\left(  1+t\right)  dt$. For
$n=0$, we obtain $J_{0}=2$ and it remains to consider the case $n\geq1$. We
employ \cite[Table 18.6.1, first row]{NIST:DLMF} to obtain%
\[
J_{n}=\int_{-1}^{1}\left(  P_{n}^{\left(  0,2\right)  }\left(  -t\right)
\right)  ^{2}\left(  1-t\right)  dt=\int_{-1}^{1}\left(  P_{n}^{\left(
2,0\right)  }\left(  t\right)  \right)  ^{2}\left(  1-t\right)  dt.
\]
From \cite[18.5.7]{NIST:DLMF} we deduce%

\begin{equation}
P_{n}^{(2,0)}\left(  t\right)  =\binom{n+2}{2}+\left(  1-t\right)
p_{n-1}\left(  t\right)  \label{expansionJacobi}%
\end{equation}
for some $p_{n-1}\in\mathbb{P}_{n-1}\left(  \left[  -1,1\right]  \right)  $.
The orthogonality property of $P_{n}^{\left(  2,0\right)  }$ leads to%
\begin{align*}
J_{n}  &  =\binom{n+2}{2}\int_{-1}^{1}P_{n}^{\left(  2,0\right)  }\left(
t\right)  \left(  1-t\right)  dt=\left(  -1\right)  ^{n}\binom{n+2}{2}%
\int_{-1}^{1}P_{n}^{\left(  0,2\right)  }\left(  -t\right)  \left(
1-t\right)  dt\\
&  =\left(  -1\right)  ^{n}\binom{n+2}{2}\int_{-1}^{1}P_{n}^{\left(
0,2\right)  }\left(  t\right)  \left(  1+t\right)  dt\\
&  =\left(  -1\right)  ^{n}\binom{n+2}{2}I_{n}=2.\hspace*{\fill}\qquad
\qquad\qquad\qquad\qquad\qquad\qquad\qquad\qquad\qquad\qquad\qquad\text{%
%TCIMACRO{\TeXButton{End Proof}{\endproof}}%
%BeginExpansion
\endproof
%EndExpansion
}%
\end{align*}

\begin{lemma}
\label{LemJacComb}For all $n\in\mathbb{N}_{0}$, it holds%
\begin{align}
\int_{-1}^{1}\left(  1-t\right)  P_{n}^{\left(  0,2\right)  }\left(  t\right)
P_{n}^{\left(  1,1\right)  }\left(  t\right)  dt  &  =2\left(  n+1\right)
,\label{0211}\\
\int_{-1}^{1}\left(  1-t\right)  P_{n}^{\left(  2,0\right)  }\left(  t\right)
P_{n}^{\left(  1,1\right)  }\left(  t\right)  dt  &  =\frac{4}{n+2}%
,\label{211b}\\
\int_{-1}^{1}P_{n}^{\left(  0,2\right)  }\left(  t\right)  dt  &  =2\left(
-1\right)  ^{n}. \label{02o}%
\end{align}

\end{lemma}

\textbf{Proof of (\ref{0211}).}

We set $T_{n}:=\int_{-1}^{1}\left(  1-t\right)  P_{n}^{\left(  0,2\right)
}\left(  t\right)  P_{n}^{\left(  1,1\right)  }\left(  t\right)  dt$. For
$n=0$, we explicitly compute $T_{0}=2$ and assume for the following $n\geq1$.
From \cite[Table 18.6.1, first row]{NIST:DLMF}, we deduce%
\[
T_{n}=\int_{-1}^{1}\left(  1+t\right)  P_{n}^{\left(  2,0\right)  }\left(
t\right)  P_{n}^{\left(  1,1\right)  }\left(  t\right)  dt.
\]
We apply expansion (\ref{expansionJacobi}) to $P_{n}^{\left(  2,0\right)
}\left(  t\right)  $ and obtain%
\[
P_{n}^{(2,0)}\left(  t\right)  =\binom{n+2}{2}+\left(  1-t\right)
p_{n-1}\left(  t\right)  \quad\text{for }p_{n-1}\in\mathbb{P}_{n-1}\left(
\left[  -1,1\right]  \right)  .
\]
The orthogonality relations of the Jacobi polynomials yield%
\begin{align}
T_{n}  &  =\binom{n+2}{2}\int_{-1}^{1}\left(  1+t\right)  P_{n}^{\left(
1,1\right)  }\left(  t\right)  dt=\left(  -1\right)  ^{n}\binom{n+2}{2}%
\int_{-1}^{1}\left(  1+t\right)  P_{n}^{\left(  1,1\right)  }\left(
-t\right)  dt\label{Tkrep}\\
&  =\left(  -1\right)  ^{n}\binom{n+2}{2}L_{n}\quad\text{with\quad}L_{n}%
:=\int_{-1}^{1}\left(  1-t\right)  P_{n}^{\left(  1,1\right)  }\left(
t\right)  dt.\nonumber
\end{align}

We compute $L_{n}$. For $n=0$ we get $L_{0}=2$ and consider $n\geq1$ in the
following. We employ \cite[18.9.5 with $\alpha=1$ and $\beta=0$ therein]%
{NIST:DLMF} to obtain%

\[
2\left(  n+1\right)  P_{n}^{(1,0)}=(n+2)P_{n}^{(1,1)}+(n+1)P_{n-1}^{(1,1)}.
\]
The orthogonality relations for Jacobi polynomials imply%
\begin{equation}
L_{n}:=-\frac{n+1}{n+2}L_{n-1}=2\frac{\left(  -1\right)  ^{n}}{n+2}%
L_{0}=\left(  -1\right)  ^{n}\frac{4}{n+2}. \label{Lk}%
\end{equation}
The combination of (\ref{Tkrep}) and (\ref{Lk}) leads to the assertion.

\textbf{Proof of (\ref{211b}).}

We set $G_{n}:=\int_{-1}^{1}\left(  1-t\right)  P_{n}^{\left(  2,0\right)
}\left(  t\right)  P_{n}^{\left(  1,1\right)  }\left(  t\right)  dt$. For
$n=0$, we obtain $G_{0}=2$ and consider $n\geq1$ in the following. From
\cite[18.5.7]{NIST:DLMF} we deduce%

\[
P_{n}^{(1,1)}\left(  t\right)  =n+1+\left(  1-t\right)  p_{n-1}\left(
t\right)  \quad\text{for some }p_{n-1}\in\mathbb{P}_{n-1}\left(  \left[
-1,1\right]  \right)  \text{.}%
\]

The orthogonality relation of $P_{n}^{\left(  2,0\right)  }$ imply%
\begin{align*}
G_{n}  &  =\left(  n+1\right)  \int_{-1}^{1}\left(  1-t\right)  P_{n}^{\left(
2,0\right)  }\left(  t\right)  dt=\left(  n+1\right)  \int_{-1}^{1}\left(
1+t\right)  P_{n}^{\left(  2,0\right)  }\left(  -t\right)  dt\\
&  =\left(  n+1\right)  \left(  -1\right)  ^{n}\int_{-1}^{1}\left(
1+t\right)  P_{n}^{\left(  0,2\right)  }\left(  t\right)
dt\overset{\text{(\ref{formel1IntP1})}}{=}\frac{4}{n+2}.
\end{align*}

\textbf{Proof of (\ref{02o}).}

Let $S_{n}:=\int_{-1}^{1}P_{n}^{\left(  0,2\right)  }\left(  t\right)  dt$.
For $n=0$, we obtain $S_{0}=2$ and assume $n\geq1$ for the following. We
employ (\ref{recJacobi}) and obtain%
\begin{equation}
S_{n}:=\int_{-1}^{1}P_{n}^{\left(  0,2\right)  }\left(  t\right)  dt=-\frac
{n}{n+2}S_{n-1}+2\frac{n+1}{n+2}H_{n}\quad\text{for }H_{n}:=\int_{-1}^{1}%
P_{n}^{\left(  0,1\right)  }\left(  t\right)  dt. \label{defSn}%
\end{equation}
We compute $H_{n}$. For $n=0$, it holds $H_{0}=2$ while for $n\geq1$ we apply
\cite[18.9.5 with $\alpha=\beta=0$.]{NIST:DLMF}, i.e.,%

\[
(2n+1)P_{n}^{(0,0)}=(n+1)P_{n}^{(0,1)}+nP_{n-1}^{(0,1)}%
\]
to obtain the recursion (using the orthogonality relation for $P_{n}^{\left(
0,0\right)  }$)%
\[
H_{n}:=-\frac{n}{n+1}H_{n-1}.
\]
Hence, $H_{n}=\frac{\left(  -1\right)  ^{n}}{n+1}H_{0}=2\frac{\left(
-1\right)  ^{n}}{n+1}$ and the recursion in (\ref{defSn}) takes the form%
\[
S_{n}=-\frac{n}{n+2}S_{n-1}+4\frac{\left(  -1\right)  ^{n}}{n+2}.
\]
This recursion with $S_{0}=2$ is satisfied for $S_{n}=2\left(  -1\right)
^{n}.$%
%TCIMACRO{\TeXButton{End Proof}{\endproof}}%
%BeginExpansion
\endproof
%EndExpansion

\bibliographystyle{abbrv}
\bibliography{nlailu}

\end{document}